\documentstyle[amstex,amssymb,12pt]{amsart}
\newcommand{\cplx}{{\Bbb C}}
\newcommand{\cz}{\cplx[z_1^{\pm 1},\dots,z_N^{\pm 1}]}
\newcommand{\zint}{{\Bbb Z}}
\newcommand{\zintp}{{\Bbb Z}_{\geq 0}}
\newcommand{\cplxn}{\cplx^{n}}
\newcommand{\ep}{\epsilon}
\newcommand{\Hhat}{\widehat{H}_N(q)}
\newcommand{\HH}{H_N(q)}
\newcommand{\sll}{\widehat{ \frak{s}\frak{l}}_n}
\newcommand{\End}{{\mathrm E}{\mathrm n}{\mathrm d}}
\newcommand{\Ker}{{\mathrm K}{\mathrm e}{\mathrm r}}
\newcommand{\sltor}{U_q'( {\frak{s}\frak{l}}_{n,tor})}
\newcommand{\Usln}{U_q'(\widehat{ \frak{s}\frak{l}}_{n})}
\newcommand{\sln}{ {{\frak s} {\frak l}}_n}
\newcommand{\glN}{ {{\frak g} {\frak l}}_N}
\newcommand{\slt}{\widehat{ \frak{s}\frak{l}}_2}

\newcommand{\mm}{{\bold m}}
\newcommand{\nn}{{\bold n}}
\newcommand{\zz}{{\bold z}}
\newcommand{\ee}{{\bold e}}
\newcommand{\ff}{{\bold f}}
\newcommand{\kk}{{\bold k}}
\newcommand{\hh}{{\bold h}}

\newcommand{\oro}{\overline{\rho}^{M,k}_{l}}
\newcommand{\orom}{\overline{\rho}^{M,k}_{l,m}}
\newcommand{\orol}{\overline{\rho}^{M,k}_{l,l+1}}

\newcommand{\MC}{{\cal M}}
\newcommand{\EC}{{\cal E}}

\newcommand{\ba}{\begin{array}}
\newcommand{\ea}{\end{array}}
\newcommand{\beq}{\begin{equation}}
\newcommand{\eeq}{\end{equation}}
\newcommand{\bqa}{\begin{eqnarray}}
\newcommand{\eqa}{\end{eqnarray}}
\newcommand{\bqas}{\begin{eqnarray*}}
\newcommand{\eqas}{\end{eqnarray*}}
\newcommand{\halmos}{\rule{5pt}{5pt}}

\newcommand{\tor}{U_q({\frak {sl}}_{n,tor})}
\newcommand{\tora}{{U_q^{*}}'({\frak {sl}}_{n,tor})}

\newcommand{\torrr}{U_q'({\frak {sl}}_{n,tor})}
\newcommand{\affa}{U_q^{(1)}(\widehat{\frak {sl}}_{n})}
\newcommand{\affaa}{U_q^{(1)'}(\widehat{\frak {sl}}_{n})}
\newcommand{\affb}{U_q^{(2)}(\widehat{\frak {sl}}_{n})}
\newcommand{\affbb}{U_q^{(2)'}(\widehat{\frak {sl}}_{n})}
\newcommand{\wb}{\widehat{B}}

\numberwithin{equation}{section}
\newtheorem{df}{\bf Definition}
\newtheorem{prop}{\bf Proposition}
\newtheorem{thm}[prop]{\bf Theorem}
\newtheorem{lemma}{\bf Lemma}
\newtheorem{cor}{\bf Corollary}
\newtheorem{con}{\bf Conjecture}
\newenvironment{rmk}{\noindent{\bf Remark}\hskip 5pt}{}

\renewenvironment{pf}{\noindent {\em {\normalsize P}\normalsize{roof.}}  
\normalsize\hskip 5pt}{\hfill\halmos}

\begin{document}
\begin{flushright} RIMS-1133 \\ \end{flushright}  \vspace{0.5cm}
\title[Toroidal actions  ]
 {Toroidal actions on level 1 modules of $U_q(\sll).$ }
\author{Yoshihisa Saito, Kouichi Takemura and Denis Uglov }
\address[ ]{Research Institute for Mathematical Sciences,
Kyoto University, Japan   }
\email{ yosihisa@@kurims.kyoto-u.ac.jp}  
\thanks{Y.S. is supported by the JSPS Research Fellowship for Young Scientist.}
\maketitle

\begin{abstract}
Recently Varagnolo and Vasserot established that the $q$-deformed Fock spaces  
due to Hayashi, and Kashiwara, Miwa and Stern, admit actions of the  quantum toroidal algebra $\sltor$ $(n\geq 3)$ with the level $(0,1)$. In the present article we propose a more detailed proof of this fact then the one given by Varagnolo and Vasserot. The proof is  based on certain non-trivial properties of  Cherednik's commuting difference operators.
The quantum toroidal action on the Fock space 
depends on a certain parameter $\kappa$. We find that 
with a specific choice of this parameter the action on the Fock spaces 
gives rise to the toroidal action on irreducible level-1 
highest weight modules of the affine quantum algebra $U_q(\sll)$. 
Similarly, by a specific choice of the parameter, the level $(1,0)$ vertex 
representation of the quantum toroidal algebra gives rise to a $\sltor$-module
structure on irreducible level-1 highest weight $U_q(\sll)$-modules.
\end{abstract}
\section{Introduction}
Recently a new algebraic object -- the quantum toroidal algebra $\sltor$ has been introduced in \cite{gkv}, \cite{VV1}. This quantum algebra is  a $q$-deformation of the enveloping algebra of a central extension of the Lie algebra $\sln[s^{\pm 1},t^{\pm 1}].$ Several results concerning representations of the quantum toroidal algebra were obtained \cite{VV1}, \cite{saito}, \cite{VV2}. One of of these results is the  Schur-type  duality between representations of the toroidal Hecke algebra and representations of $\sltor$ established by  Varagnolo and Vasserot \cite{VV1}. This duality is analogous to the duality between affine Hecke algebra and the quantum affine algebra $U_q'(\sll)$ given by Chari and Pressley \cite{CP}. It is known that $\sltor$ contains two subalgebras $\affaa$ and $\affbb$ isomorphic to $\Usln.$ A module of $\sltor$ is said to have level $(k,l)$ if the $\affaa$-action on this module has level $k$ and the  $\affbb$-action has level $l.$ 

Representations of $\sltor$ obtained by  the Varangolo-Vasserot duality have the level $(0,0).$ They are analogues of level-0 modules of affine quantum algebras.     

The first example of a toroidal module with a non-trivial level was constructed in \cite{saito}. This vertex operator construction is   an analogue of the Frenkel-Jing bosonization for level-1 $U_q'(\sll)$-modules. It gives a toroidal module with level $(1,0).$ We summarize this bosonic construction in section \ref{sec:saito}.         

A $q$-fermionic construction of toroidal modules with non-trivial level has recently been proposed in \cite{VV2}. Origins of this construction lay in the theory of integrable models with long-range interaction. It has been known for some time, starting with the work \cite{BPS}, that level-1 highest weight modules of  $U_q'(\sll)$ admit a level-0 action of the same quantum affine algebra $U_q'(\sll)$ \cite{JKKMP},\cite{TU}. In particular in   \cite{TU} it was shown that the Fock space module of  $U_q'(\sll)$ \cite{H, KMS} is simultaneously a level-0  $U_q'(\sll)$-module.          

The main result of the work \cite{VV2} is that for $n \geq 3$ the level-0 and level-1 actions of $U_q'(\sll)$ on the Fock space are exactly actions of the subalgebras  
$\affaa$ and $\affbb$ in the toroidal algebra, so that the Fock space may be regarded as a level $(0,1)$ module of $\sltor.$ 

\mbox{}From our viewpoint the proof of this fact given in  the Theorem  of  section 12 in  the paper \cite{VV2} omits certain technical details -- notably in the  proof of the relations (\ref{E1}) - (\ref{Kn-1}) of the present work. One of our main objectives in this paper is to supply these details by giving a  complete  proof, based on Lemmas \ref{l: lemma3} and \ref{stab},   that the Fock space is a module of the toroidal quantum algebra (Theorem \ref{mthm}). Let us emphasize that the general idea of our proof is not new methodologically compared with the idea of the proof in \cite{VV2}. We arrived at this idea before the appearance of the work \cite{VV2}. However, at the time when the paper \cite{VV2} became available to us we did not know a complete proof of Theorem \ref{mthm}, specifically a proof of the relations (\ref{En-1}) - (\ref{Kn-1}) in Lemma \ref{stab} was missing. 

In this paper we will also show that  the vertex representation  of \cite{saito} is isomorphic to 
the Fock space as a level-1  $U_q'(\sll)$-module. This means that on the Fock 
space we have two actions of the toroidal algebra -- one action of level 
$(1,0)$ and another action of level $(0,1)$ such that the action of the 
$\affaa $-subalgebra in the former coincides with the action of the 
$\affbb$-subalgebra in the latter. We are not aware of any good explanation 
of this phenomenon.

Apart from the  matters discussed above  we give a proof of irreducibility of the Fock space as a $\sltor$-module. We also demonstrate that irreducible highest weight level-1 modules 
of $U_q'(\sll)$ admit actions of toroidal algebra with levels $(1,0)$ and $(0,1)$  induced from the corresponding actions on the Fock space.

\section{Definition of quantum toroidal algebras}
\subsection{ }
Let ${\frak {sl}}_n$ be the semisimple Lie algebra of type $A_{n-1}$ and
$\widehat{\frak {sl}}_n$ the affine Kac-Moody Lie algebra of type 
$A_{n-1}^{(1)}$.
We denote their  Cartan subalgebras by ${\frak h}$ and $\widehat{\frak h}$.
We denote by ${\alpha}_0,\cdots,{\alpha}_{n-1}$ the simple
roots and by ${h}_0,\cdots,{h}_{n-1}$ the simple coroots of $\widehat{ \frak{s}\frak{l}}_{n}.$ Let 
$P=\oplus_{i=0}^{n-1}{\Bbb Z}\Lambda_i\oplus{\Bbb Z}\delta$ be the weight 
lattice. Here $\delta$ is the null root. Let $Q=\oplus_{i=0}^{n-1}
{\Bbb Z}\alpha_i$ be the root lattice. Note that $\alpha_0=-\Lambda_{n-1}
+2\Lambda_0-\Lambda_1+\delta$, $\alpha_i=-\Lambda_{i-1}+2\Lambda_i
-\Lambda_{i+1}$ $(1\leq i\leq n-1)$. Here the indices are extended cyclically
such as $\Lambda_i=\Lambda_{i+n}$. Let $\overline{P}=\oplus_{i=1}^{n-1}{\Bbb Z}
\overline{\Lambda}_i$ be the classical weight lattice and $\overline{Q}=
\oplus_{i=1}^{n-1}{\Bbb Z}\alpha_i$ the classical root lattice. The inclusion of $\overline{P}$ into  $P$ is given by $\overline{\Lambda_i}=\Lambda_i
-\Lambda_0.$ We also set $\overline{\delta}=0$.

We denote the pairing of $\frak h$ and
${\frak h}^*$ ({\it resp}. $\widehat{\frak h}$ and $\widehat{\frak h}^*$) by
$\langle\text{ },\text{ }\rangle$. The invariant bilinear form on $P$ is
given by $(\alpha_i|\alpha_j)=-\delta_{ij-1}+2\delta_{ij}-\delta_{ij+1}$ and
$(\delta|\delta)=0$. 

Let $[k]=(q^k-q^{-k})/({q-q^{-1}})$ be the $q$-integers and  denote
$[n]!=\prod_{k=1}^n[k]$, $\begin{bmatrix}m\\r\end{bmatrix}
={[m]!}/({[r]![m-r]!})$.
\subsection{ }
$U_q(\widehat{\frak{sl}}_n)$ is a ${\Bbb Q}(q)$-algebra generated by the 
symbols $E_i$, $F_i$, $K_i^{\pm}$ $(i=0,\cdots,n-1)$ $q^{\pm d}$ which satisfy
the following defining relations:
$$K_i^{\pm}K_j^{\pm}=K_j^{\pm}K_i^{\pm},\quad K_i^+E_jK_i^-=
q^{\langle h_i,\alpha_j\rangle}E_j,
\quad K_i^+F_jK_i^-=q^{-\langle h_i,\alpha_j\rangle}F_j,$$
$$K_i^{\pm}q^{\pm d}=q^{\pm d}K_i^{\pm},\quad q^dE_jq^{-d}=q^{\delta_{0j}}E_j,
\quad q^dF_jq^{-d}=q^{-\delta_{0j}}F_j,$$
$$[E_i,F_j]=\delta_{ij}\frac{K_i^+-K_i^-}{q-q^{-1}},$$
for $i\ne j$,
$$\sum_{r=0}^{m}(-1)^{r}\begin{bmatrix}m\\r\end{bmatrix}E_i^rE_jE_i^{m-r},
\quad \sum_{r=0}^{m}(-1)^{r}\begin{bmatrix}m\\r\end{bmatrix}F_i^rF_jF_i^{m-r}.
$$
Here $m=\langle h_i,\alpha_j\rangle$.

Let $\Usln$ be the subalgebra generated by $E_i$, $F_i$, $K_i^{\pm}$ 
$(i=0,\cdots,n-1)$. 

The coproduct $\Delta$ of $U_q(\widehat{\frak{sl}}_n)$ is given as follows:
$$\Delta(K_i^{\pm})=K_i^{\pm}\otimes K_i^{\pm},\quad 
  \Delta(q^{\pm d})=q^{\pm d}\otimes q^{\pm d},$$
$$\Delta(E_i)=E_i\otimes K_i^++1\otimes E_i,\quad \Delta(F_i)=F_i\otimes 1+
K_i^-\otimes F_i.$$

\subsection{ }
We will give the definition of the quantum toroidal algebra $\tor$. Fix an 
integer $n\geq 3$ and $\kappa\in {\Bbb Q}(q)^{\times}$.
\begin{df}
The quantum toroidal algebra $\tor$ is an associative algebra over 
${\Bbb Q}(q)$ with generators :
$$E_{i,k},\text{ }F_{i,k},\text{ }H_{i,l},\text{ }K_i^{\pm},\text{ }
q^{\pm\frac12 c},q^{\pm d_1},q^{\pm d_2},$$
for $k\in {\Bbb Z}$, $l\in {\Bbb Z}\backslash \{0\}$ and $i=0,1,\cdots,n-1$.\\

The relations are expressed in terms of the formal series
$$E_i(z)=\sum_{k\in {\Bbb Z}}E_{i,k}z^{-k},$$
$$F_i(z)=\sum_{k\in {\Bbb Z}}F_{i,k}z^{-k},$$
$$K_i^{\pm}(z)=K_i^{\pm}\exp(\pm (q-q^{-1})\sum_{k\geq 1}H_{i,\pm k}
z^{\mp k}),$$
as follows:
\begin{equation}
 q^{\pm\frac12 c}\text{ are central,}
\label{defrel1}
\end{equation}
\begin{equation}
 K_i^+K_i^{-}=K_i^{-}K_i^+=1,
\end{equation}
\begin{equation}
 K_i^{\pm}(z)K_j^{\pm}(w)=K_j^{\pm}(w)K_i^{\pm}(z)
\label{rb}
\end{equation}
\begin{equation}
 \theta_{-\langle h_i,\alpha_j\rangle}
 (q^{-c}\kappa^{m_{ij}}\frac{z}{w})K_i^-(z)K_j^+(w)=
 \theta_{-\langle h_i,\alpha_j\rangle}
 (q^{c}\kappa^{m_{ij}}\frac{z}{w})K_j^+(w)
 K_i^-(z)
\label{rc}
\end{equation}
\begin{equation}
q^{d_1}K_i^{\pm}(z)q^{-d_1}=K_i^{\pm}(q^{-1}z),\quad
q^{d_1}E_i(z)q^{-d_1}=E_i(q^{-1}z),\quad
q^{d_1}F_i(z)q^{-d_1}=F_i(q^{-1}z),
\end{equation}
$$
[q^{d_2},K_i^{\pm}(z)]=0,\quad
q^{d_2}E_i(z)q^{-d_2}=q^{\delta_{i0}}E_i(z),\quad
q^{d_2}F_i(z)q^{-d_2}=q^{-\delta_{i0}}F_i(z),
$$
\begin{equation}
 K_i^{\pm}(z)E_j(w)
 =\theta_{\mp\langle h_i,\alpha_j\rangle}
 (q^{-\frac12 c}\kappa^{\mp m_{ij}}w^{\pm}z^{\mp})E_j(w)K_i^+(z)
\label{rd} \end{equation}
$$
 K_i^{\pm}(z)F_j(w)
 =\theta_{\pm\langle h_i,\alpha_j\rangle}
 (q^{\frac12 c}\kappa^{\mp m_{ij}}w^{\pm}z^{\mp})F_j(w)K_i^+(z)
$$
\begin{equation}
 [E_i(z),F_j(w)]=\delta_{i,j}\frac{1}{q-q^{-1}}
 \{\delta(q^c\frac{w}{z})K_i^+(q^{\frac12 c}w)-\delta(q^c\frac{z}{w})K_i^
 -(q^{\frac12 c}z)\}
\label{re}
\end{equation}
\begin{equation}
 (\kappa^{m_{ij}}z-q^{\langle h_i,\alpha_j\rangle}w)E_i(z)E_j(w)
 =(q^{\langle h_i,\alpha_j\rangle}\kappa^{m_{ij}}z-w)E_j(w)E_i(z)
\label{rf}
\end{equation}

$$
 (\kappa^{m_{ij}}z-q^{-\langle h_i,\alpha_j\rangle}w)F_i(z)F_j(w)
 =(q^{-\langle h_i,\alpha_j\rangle}\kappa^{m_{ij}}z-w)F_j(w)F_i(z)
$$
\begin{equation}
 \sum_{\sigma\in {\frak S}_m}\sum_{r=0}^m(-1)^r
 \begin{bmatrix}m\\r\end{bmatrix}
 E_i(z_{\sigma(1)})\cdots E_i(z_{\sigma(r)})E_j(w)E_i(z_{\sigma(r+1)})\cdots
 E_i(z_{\sigma(m)})=0
\label{rg}
\end{equation}
\begin{equation*}
  \sum_{\sigma\in {\frak S}_m}\sum_{r=0}^m(-1)^r
 \begin{bmatrix}m\\r\end{bmatrix}
{} F_i(z_{\sigma(1)})\cdots F_i(z_{\sigma(r)})F_j(w)F_i(z_{\sigma(r+1)})\cdots
{} F_i(z_{\sigma(m)}) = 0 
\end{equation*}

where $i\ne j$ and $m=1-\langle h_i,\alpha_j\rangle$.\\

In these formulas we denote $\theta_m(z)=\frac{zq^m-1}{z-q^m}$ for
$m\in {\Bbb Z}$, the $\theta_m(z)$ is to be regarded as the expansion of the right-hand side above in non-negative powers of the argument $z$;  $\delta(z)=\sum_{k\in {\Bbb Z}}z^k$, $m_{ij}$ are the entries
of the following $n\times n$-matrix
$$
M=\begin{pmatrix}
      0 &     -1 &      0 & \hdots &      0 &      1\\
      1 &      0 &     -1 & \hdots &      0 &      0\\
      0 &      1 &      0 & \hdots &      0 &      0\\
 \vdots & \vdots & \vdots & \ddots & \vdots & \vdots\\
      0 &      0 &      0 & \hdots &      0 &     -1\\
     -1 &      0 &      0 & \hdots &      1 &      0
\end{pmatrix}.
$$
\end{df}
Let $\torrr$ be the subalgebra of $\tor$ generated by $E_{i,k},F_{i,k},
K_i^{\pm},H_{i,l}$, $q^{\pm\frac12 c}$. 

\subsection{}
Let $\affaa$ be the subalgebra of $\tor$ generated
by $E_{i,k},F_{i,k},K_i^{\pm},H_{i,l}$, $q^{\pm\frac12 c}$
$(1\leq i \leq n-1,k\in{\Bbb Z},
l\in{\Bbb Z}\backslash\{0\})$ and $\affa$ the subalgebra generated by
$\affaa$ and $q^{\pm d_1}$.
Let $\affbb$ the subalgebra of $\tor$ generated by $E_{i,0},F_{i,0},K_i^{\pm}$
$(0\leq i \leq n-1)$ and $\affb$ the subalgebra generated by
$\affbb$ and $q^{\pm d_2}$.

The following lemma is already known \cite{gkv}, \cite{VV1}.
\begin{lemma}
 Both $\affa$ and $\affb$ are isomorphic to $U_q(\widehat{\frak {sl}}_{n})$.
\end{lemma}
Following \cite{VV2} we will call the $\affaa$ the vertical subalgebra of $\sltor$ and denote its Chevalley generators by $\ee_i , \ff_i, \kk_i^{\pm 1}$ $ (i \in \{0,\dots,n-1\}).$ The $\affbb$ will accordingly be called the horizontal subalgebra and for its Chevalley generators $E_{i,0},F_{i,0},K_i^{\pm}$ we will use the notations $ E_{i},F_{i},K_i^{\pm}$ $ (i \in \{0,\dots,n-1\})$ respectively.

\begin{lemma}  \label{l: gen}
 $\affa$ and $\affb$ generate $\tor$.
\end{lemma}
\begin{pf}
 Let $\cal U$ be the subalgebra of $\tor$ which is generated by $\affa$ and 
 $\affb$. It is enough to show that $E_{0,k},F_{0,k},H_{0,l}\in {\cal U}$.
 By the defining relations we have
 $$
  [H_{1,l},E_{0,k}]=-\frac{1}{l}[~l~]q^{-\frac12 |l|c}\kappa^{-l}E_{0,k+l}.  
 $$
 Since $E_{0,0}$ and $H_{1,l}$ $(l\in {\Bbb Z}\backslash \{0\})$ are the 
 elements of $\cal U$ we have $E_{0,k}\in {\cal U}$ for any $k$. Similarly
 we have $F_{0,k}, H_{0,l}\in {\cal U}$. 
\end{pf}

\begin{df}
 Let $V$ be a $\tor$-module. We say that $V$ has level $(l_1,l_2)$ if $V$ has level 
 $l_j$ as $U_q^{(j)}(\widehat{\frak {sl}}_{n})$-module $(j=1,2)$. 
\end{df}
On a level $(l_1,l_2)$ module the generator $q^{\frac{1}{2}c}$ acts as the multiplication by $q^{\frac{1}{2}l_2}$ and the element $ K_0 K_1 \dots K_{n-1} $ acts as the multiplication by $q^{l_1}.$

\section{The vertex representation} \label{sec:saito}
\subsection{ }
In this section we assume $c=1$. Let us recall the results on the vertex
representation of the quantum toroidal algebra~\cite{saito}.

Let $S_n$ be the subalgebra of $\tor$, generated by the $H_{i,l}$ 
$(0\leq i \leq n-1,\text{ } l\in {\Bbb Z}\backslash \{0\})$. By definition
the commutation relations of $\{ H_{i,l}\}$ are following:
$$
[H_{i,k},H_{j,l}]=\delta_{k+l,0}\frac1k [k\langle h_i,\alpha_j\rangle]
\frac{q^{k}-q^{-k}}{q-q^{-1}}\kappa^{-km_{ij}}. 
$$
Let ${\cal F}_n$ be the Fock space of $S_n$. That is, ${\cal F}_n$ is
generated by the vacuum vector $v_0$ and the defining relations are 
$H_{i,l}v_0=0$ for $l>0$.
\subsection{ }   
We introduce a twisted  group algebra ${\Bbb Q}(q)\{\overline{P}\}$ defined as 
the ${\Bbb Q}(q)$-algebra generated by symbols $e^{\pm {\alpha_2}}$,
$e^{\pm {\alpha_3}},
\cdots,\text e^{\pm {\alpha_{n-1}}}$, $e^{\pm \overline{\Lambda}_{n-1}}$
which satisfy the following relations:
\begin{equation*}
 e^{{\alpha_i}}e^{{\alpha_j}}=(-1)^{\langle {h_i},
 {\alpha_j}\rangle}e^{{\alpha_j}}e^{{\alpha_i}}, \qquad 
 e^{{\alpha_i}}e^{\overline{\Lambda}_{n-1}}=(-1)^{\delta_{i,n-1}}
  e^{\overline{\Lambda}_{n-1}}e^{{\alpha_i}}.
\end{equation*}
For ${\alpha}=\sum_{i=2}^{n-1} m_i {\alpha_i}+m_{n}
\overline{\Lambda}_{n-1}$ we denote
$e^{{\alpha}}=(e^{{\alpha_2}})^{m_2}
(e^{{\alpha_3}})^{m_3}
\cdots (e^{{\alpha}_{n-1}})^{m_{n-1}}
(e^{\overline{\Lambda}_{n-1}})^{m_{n}}$.

We denote by ${\Bbb Q}(q)\{\overline{Q}\}$ the subalgebra of
${\Bbb Q}(q)\{\overline{P}\}$ generated by $e^{{\alpha_i}}$
($1\leq i \leq n-1$). \\
Set
$$W(M)={\cal F}_n\otimes {\Bbb Q}(q)\{\overline{Q}\}
e^{\overline{\Lambda}_M}\qquad\text{for }0\leq M\leq n-1,$$
here we denote $\overline{\Lambda}_0=0$.

For $i=0,1,\cdots,n-1$ we denote 
$$
\overline{\alpha}_i=
\begin{cases}
 -\sum_{j=1}^{n-1}{\alpha}_j, & \quad i=0,\\
 \alpha_i, & \quad i\ne 0,
\end{cases}
\qquad
\overline{h}_i=
\begin{cases}
 -\sum_{j=1}^{n-1}h_j, & \quad i=0,\\
 h_i, & \quad i\ne 0,
\end{cases}
$$

We define the operators on $W(M)$, $H_{i,l}$, $e^{{\alpha}}$, 
$({\alpha}\in \overline{Q})$, $\partial_{\overline{\alpha}_i}$, $z^{H_{i,0}}$ 
$(i=0,1,\cdots,n-1)$ and $d$
as follows:\\
For $v\otimes e^{\beta}=H_{i_1,-k_1}\cdots H_{i_N,-k_N}v_0
\otimes e^{\sum_{j=1}^{n-1}m_j\alpha_j+\overline{\Lambda}_M}\in W(M)$,
$$H_{i,l}(v\otimes e^{{\beta}})=(H_{i,l}v)\otimes
e^{{\beta}},$$
$$e^{\alpha}(v\otimes e^{\beta})=v\otimes e^{\alpha}e^{\beta},$$
$${\partial_{\overline{\alpha_i}}}(v\otimes e^{\beta})=
 {\langle\overline{h_i},{\beta}\rangle}v\otimes e^{\beta},$$
$$z^{H_{i,0}}(v\otimes e^{\beta})
 =z^{\langle \overline{h_i},\beta\rangle}
 \kappa^{\frac12 \sum_{k=1}^{n-1}\langle \overline{h_i},m_k\alpha_k\rangle 
 m_{ik}}(v\otimes e^{\beta}),$$
$$d(v\otimes e^{{\beta}})=(-\sum_{s=1}^{N}k_s
 -\frac{({\beta}|{\beta})}{2}
 +\frac{(\overline{\Lambda}_M|\overline{\Lambda}_M)}{2})v\otimes e^{\beta}.$$
Let us denote by $\tora$ 
the subalgebra of $\tor$ generated by $E_{i,k},F_{i,k},K_i^{\pm},H_{i,l}$, 
$q^{\pm\frac12 c}$ and $q^{\pm d_1}.$ 
\begin{prop} \cite{saito}
 Let $c=1$. Then
 for each $M$ and $n$, the following action gives a $\tora$-module structure
 on $W(M)$:
 $$q^{\frac12 c}\mapsto q^{\frac12},$$
 $$q^{d_1}\mapsto q^{d},$$
 $$E_i(z)\mapsto\exp(\sum_{k\geq 1}\frac{H_{i,-k}}{[k]}(q^{-1/2}z)^k)
         \exp(\sum_{k\geq 1}-\frac{H_{i,k}}{[k]}(q^{1/2}z)^{-k})
         e^{\overline{\alpha_i}}z^{\partial_{\overline{\alpha_i}}+1},$$
 $$F_i(z)\mapsto\exp(\sum_{k\geq 1}-\frac{H_{i,-k}}{[k]}(q^{1/2}z)^k)
         \exp(\sum_{k\geq 1}\frac{H_{i,k}}{[k]}(q^{-1/2}z)^{-k})
         e^{-\overline{\alpha_i}}z^{-\partial_{\overline{\alpha_i}}+1},$$
 $$K_i^+(z)\mapsto \exp((q-q^{-1})
   \sum_{k\geq 1}H_{i,k}z^{-k})q^{\partial_{\overline{\alpha_i}}},$$
 $$K_i^-(z)\mapsto \exp(-(q-q^{-1})
   \sum_{k\geq 1}H_{i,-k}z^{k})q^{-\partial_{\overline{\alpha_i}}}$$
 for $0\leq i\leq n-1$. Moreover $W(M)$ has level $(1,0)$ as a $\tora$-module. 
\end{prop}

\begin{prop} \cite{saito} \label{prop;ch}
 As a $\affa$-module 
 $$\text{ch}_{W(M)}=\frac{\text{ch}_{L(\Lambda_M)}}{\varphi(e^{-\delta})},$$
 where $L(\Lambda_M)$ is the irreducible highest weight 
 $U_q(\widehat{\frak{sl}_n})$-module with the highest weight $\Lambda_M$ and
 $\varphi(x)=\prod_{k>0}(1-x^k)$.  
\end{prop}

\begin{df}
 We say that $\kappa\in {\Bbb Q}(q)^{\times}$ is generic if for any 
 $k\in{\Bbb Z}_{>0}$
 the $n\times n$-matrix $G(n,k,\kappa)=([k\langle h_i,\alpha_j\rangle]
 \kappa^{-km_{ij}})$ is invertible.
\end{df}

\begin{thm} \cite{saito}
 If $\kappa$ is generic then $W(M)$ is irreducible $\tora$-module.
\end{thm}
\subsection{Non-generic case}
We shall start the following lemma.
\begin{lemma}
 If $\kappa$ is not generic then $\kappa \in \{\pm q, \pm q^{-1}\}.$
\end{lemma}
\begin{pf}
 We assume $\text{det}G(n,k,\kappa)=0$ for fixed $k$. By an easy calculation 
 we have $\text{det}G(n,k,\kappa)=q^{-nk}[k](q^{nk}-\kappa^{nk})
 (q^{nk}-\kappa^{-nk})$. 
 Therefore if $nk$ is even then $\kappa=\pm q^{\pm 1}$ and if $nk$ is odd then
 $\kappa =q^{\pm 1}$. Thus we have the statement. 
\end{pf}

Until the end of this subsection we assume $\kappa=q^{\pm 1}$.
Let $\wb_k=\sum_{i=0}^{n-1}H_{i,k}$ and $\widehat{S}$ be the algebra generated
by $\wb_k$ $(k\in {\Bbb Z}\backslash \{0\})$. It is easy to see the following 
lemma. 
\begin{lemma}
 If $\kappa=q^{\pm 1}$ then $[H_{i,l},\wb_k]=0$ for any $i,k$ and $l$. 
 Moreover $\widehat{S}$ is an abelian subalgebra of $\tor$. 
\end{lemma}

By the above lemma and the definition of the action on $W(M)$ we have the 
following lemma immediately.
\begin{lemma}\label{lemma;comm}
 $[\widehat{S},\torrr]=0$ on $W(M)$. 
\end{lemma}

Let $\widehat{S}^{<0}=\oplus_{k< 0} {\Bbb Q}(q) \{X\in \widehat{S}|q^{d_1}Xq^{-d_1}=q^{k}X \}$. Set $W(M)^{(1)}=\widehat{S}^{<0}W(M)$.
\begin{prop}
 If $\kappa=q^{\pm 1}$ then \\
 $(1)$ $W(M)^{(1)}$ is a proper $\tora$-submodule of $W(M)$ and\\
 $(2)$ $W(M)/W(M)^{(1)}$ is an irreducible $\tora$-module.
\end{prop}
\begin{pf}
 $(1)$ By Lemma \ref{lemma;comm} $W(M)^{(1)}$ has a $\tora$-module structure. 
 Since $v_0\otimes e^{\overline{\Lambda}_M}$ does  not belong to $W(M)^{(1)}$, 
 it is a proper submodule.\\
 $(2)$ As $\affa$-module the character of $W(M)$ is already known 
 (See Proposition \ref{prop;ch}). By the definition we have
 $$\text{ch}_{W(M)/W(M)^{(1)}}=\text{ch}_{L(\Lambda_M)}$$
 as $\affa$-module. Thus ${W(M)/W(M)^{(1)}}\cong L(\Lambda_M)$
 as $\affa$-module. Therefore $W(M)/W(M)^{(1)}$ is irreducible. 
\end{pf}
\begin{cor}
$(1)$ If $\kappa=q^{\pm 1}$ then $L(\Lambda_M)$ has $\tor$-module structure.\\
$(2)$ Therefore we have another $U_q(\widehat{\frak{sl}_n})$-module structure
on $L(\Lambda_M)$ with level $0$.
\end{cor}

\section{Action of $\sltor$ on the finite wedge product  }
\subsection{Toroidal Hecke algebra}
Let $q \in \cplx^{\times}$. The toroidal Hecke algebra of type $\glN$, $H_{tor}$ is a unital associative algebra over $\cplx[x^{\pm 1},y^{\pm 1}]$ with generators $T_i^{\pm 1},$ $X_j^{\pm 1},$  $Y_j^{\pm 1},$ $i=1,2,\dots,N-1,$ $j=1,2,\dots,N$ and relations 
\begin{gather*}
T_i T_i^{-1} =  T_i^{-1} T_i = 1, \quad (T_i + 1) (T_i - q^2) =0, \\
T_i T_{i+1} T_i =  T_{i+1} T_i T_{i+1},\\
T_i T_j = T_j T_i \quad \text{ if $ |j-i| > 1,$} \\
X_0 Y_1 = x Y_1 X_0, \quad X_i X_j = X_j X_i, \quad Y_i Y_j = Y_j Y_i, \\ 
X_j T_i = T_i X_j, \quad Y_j T_i = T_i Y_j \quad \text{ if $ j \neq i,i+1 > 1,$} \\
T_i X_i T_i = q^2 X_{i+1}, \quad T_i^{-1} Y_i T_i^{-1} = q^{-2} Y_{i+1}, \\
X_2 Y_1^{-1} X_2^{-1} Y_1 = q^{-2} y T_1^2,
\end{gather*}
where $X_0 = X_1 X_2 \cdots X_N.$ In $H_{tor}$ the subalgebras $\Hhat^{(1)}$ generated by $ T_i^{\pm 1},$ $Y_j^{\pm 1},$ and $\Hhat^{(2)}$ generated by $ T_i^{\pm 1},$ $X_j^{\pm 1}$ are both isomorphic to the affine Hecke algebra. And the subalgebra $\HH$ generated by $ T_i^{\pm 1}$ is isomorphic to the Hecke algebra of type $\glN.$

Let $p \in \cplx^{\times}$ and consider the following operators in $\End(\cz)$ 
\begin{eqnarray}
& \text{multiplication operator $z_i,$}   \;\; & i=1,2,\dots,N, \nonumber\\
&g_{i,j} = \frac{q^{-1}z_i - q z_j}{z_i - z_j}(K_{i,j} - 1) + q, \;\; & 1\leq i \neq j \leq N,\nonumber \\ 
& \text{ the family of $N$ commutative Cherednik's operators: } \nonumber \\
& Y_i^{(N)} = g_{i,i+1}^{-1}K_{i,i+1} \cdots g_{i,N}^{-1} K_{i,N} p^{D_i} K_{1,i} g_{1,i} \cdots K_{i-1,i}g_{i-1,i}, \;\; & i=1,2,\dots,N, \nonumber
\end{eqnarray}
where $K_{i,j}$ acts on $\cz$ by permuting variables $z_i,z_j$ and $p^{D_i}$ is the difference operator  
\begin{equation}
p^{D_i}f(z_1,\dots,z_i,\dots,z_N)  =  f(z_1,\dots,pz_i,\dots,z_N) , \quad f \in \cplx[z_1^{\pm 1},\dots, z_N^{\pm 1}]. \nonumber 
\end{equation}
The following result is due to I. Cherednik \cite{Cherednik,Cherednik2}:
\begin{prop}  \label{p: heckeact}
The map 
\begin{equation}
T_i \mapsto \overset{c}{T_i} = -q g_{i,i+1}^{-1}, \quad X_i \mapsto z_i, \quad Y_i \mapsto q^{1-N}Y_i^{(N)}, \quad x \mapsto p, \quad y \mapsto 1  
\end{equation}
defines a right action of $H_{tor}$  on $\cz.$
\end{prop}
The commuting difference operators $Y_1^{(N)}, \dots, Y_N^{(N)}$ are called Cherednik's operators. 

Let $V = \cplxn$, with basis $\{ v_1,\dots,v_n \}.$ Then $\otimes^N V$ admits a left $\HH$-action given by    
\begin{gather}
T_i \mapsto \overset{s}{T_i} = 1^{\otimes^{i-1}}\otimes \overset{s}{T} \otimes 1^{\otimes^{N-i-1}}, \quad \text{ where } \quad  \overset{s}{T}  \in \End(\otimes^2 V) \\
\text{ and } \quad   \overset{s}{T}(v_{\ep_1}\otimes v_{\ep_2}) = \begin{cases} q^2 v_{\ep_1}\otimes v_{\ep_2} & \text{ if $ \ep_1 = \ep_2,$} \\
 q v_{\ep_2}\otimes v_{\ep_1} & \text{ if $ \ep_1 < \ep_2,$ }\\
q v_{\ep_2}\otimes v_{\ep_1} + (q^2  - 1) v_{\ep_1}\otimes v_{\ep_2}& \text{ if $ \ep_1 > \ep_2.$ }  \end{cases}
\end{gather}

\subsection{$q$-wedge product} \label{sec: qwedge}

Let $V(z)= \cplx [z^{\pm 1}]\otimes V,$ with  basis $ \{ z^m \otimes v_{\ep} \}$, $ m\in \zint $ , $ \ep \in \{1,2,\dots,n\} $. Often it will be convenient to set $ k = \ep - nm $ and $ u_k = z^m \otimes v_{\ep} $. Then $\{ u_k \}, $  $k \in \zint $ is a basis of $V(z)$. In what follows we will write $ z^m v_{\ep} $ as a short-hand for $ z^m \otimes v_{\ep} ,$ and use both  notations: $ u_k $ and $ z^m  v_{\ep} $ switching between them according to convenience. 

The two actions of the Hecke algebra are naturally extended on the tensor product $\cz\otimes(\otimes^N V)$ so that $\overset{c}T_i$ acts trivially on $\otimes^N V$ and  $\overset{s}T_i$ acts trivially on $\cz.$ The vector space $\otimes^N V(z)$ is identified with  $\cz\otimes(\otimes^N V)$ and the $q$-wedge product \cite{KMS} is defined as the following quotient space:
\begin{equation}
\wedge^N V(z) = \otimes^N V(z)/ \sum_{i=1}^{N-1}\Ker\left( \overset{c}T_i + q^2 (\overset{s}T_i)^{-1} \right). \label{e: q-wedgeprod}
\end{equation}
Since for any $i =1,2,\dots,N-1$ we have 
\begin{equation}
\Ker\left( \overset{c}T_i + q^2 (\overset{s}T_i)^{-1} \right) = {\mathrm I}{\mathrm m}\left( \overset{c}T_i - \overset{s}T_i \right)
\end{equation}
the definition (\ref{e: q-wedgeprod}) is equivalent to 
\begin{equation}
\wedge^N V(z) = \cz \otimes_{\HH} (\otimes^N V). 
\end{equation}

Let $\Lambda: \otimes^N V(z) \rightarrow \wedge^N V(z) $ be the quotient map specified by (\ref{e: q-wedgeprod}). The image of a pure tensor $ u_{k_1}\otimes u_{k_2}\otimes \dots \otimes u_{k_N} $ under this map is called a wedge and is denoted by
\begin{equation} 
u_{k_1}\wedge u_{k_2}\wedge \dots \wedge u_{k_N} :=\Lambda ( u_{k_1}\otimes u_{k_2}\otimes \dots \otimes u_{k_N}) . \label{e: wedge} 
\end{equation}
A wedge is normally ordered if $k_1 > k_2 > \cdots > k_N.$ In \cite{KMS} it is proven that normally ordered wedges form a basis in $\wedge^N V(z),$ and that any wedge is expressed as a linear combination of normally ordered wedges by using the normal ordering rules:   
\bqa
u_l\wedge u_m & = & -u_m\wedge u_l, \quad \text{for $ l=m\bmod n $}, \label{e: no1}\\
u_l\wedge u_m & = & -qu_m\wedge u_l + (q^2 - 1)(u_{m-i}\wedge u_{l+i} - qu_{m-n}\wedge u_{l+n} + \nonumber\\
& & \qquad \qquad + q^2 u_{m-n-i}\wedge u_{l+n+i} + \dots ), \label{e: no2}\\ 
& & \text{for $ l < m , m - l =i \bmod n , 0<i<n$}. \nonumber   
\eqa
The sum above continues as long as the wedges in the right-hand side are normally ordered.

\subsection{Action of the quantum toroidal algebra on the wedge  product} \label{sec: finqta}

In the paper \cite{VV1} the following result is proven
\begin{thm}[Varagnolo and Vasserot] \label{t: VV}
Suppose that $x= \kappa^{-n}q^n$ and $y=1.$ Then for any right $H_{tor}$-module $R$ the vector space $ R\otimes_{\HH} (\otimes^N V) $ is a left $\sltor$-module with  central charge $(0,0).$
\end{thm}
Moreover the action of the vertical subalgebra $\affaa $  generated by $ \ee_i, \ff_i, \kk_i^{\pm 1} $ and the action of the horizontal subalgebra  $\affbb $ generated by $ E_i, F_i, K_i^{\pm 1} $ $( i =0,1,\dots,n-1)$  are  given on the space $ R\otimes_{\HH} (\otimes^N V) $ as follows: Let $m \in R$ and $v \in \otimes^N V,$ then 
\begin{align}
& E_i ( m \otimes v)  =  \ee_i ( m \otimes v)   = \sum_{j=1}^N m \otimes E_j^{i,i+1} K_{j+1}^{i}\dots K_N^{i} v,  \label{e: Efin}\\  
& F_i ( m \otimes v)    =  \ff_i ( m \otimes v)   = \sum_{j=1}^N m \otimes (K_{1}^{i})^{-1}\dots (K_{j-1}^{i})^{-1}E_j^{i+1,i} v,  \label{e: Ffin}\\  
& K_i( m \otimes v)     =  \kk_i ( m \otimes v)    = m \otimes   K^{i}_1K^{i}_2\dots K^{i}_N v ,\quad    (i=1,2,\dots,n-1)  \label{e: Kfin} \\ 
& \ee_0( m \otimes v) =  \sum_{j=1}^N m Y_j^{-1}\otimes E_j^{n,1} K_{j+1}^{0}\dots K_N^{0} v,  \label{e: E0} \\ 
& \ff_0 ( m \otimes v)    =   \sum_{j=1}^N m Y_j \otimes (K_{1}^{0})^{-1}\dots (K_{j-1}^{0})^{-1}E_j^{1,n} v,  \label{e: F0} \\  
& E_0( m \otimes v) =  \sum_{j=1}^N m X_j \otimes E_j^{n,1} K_{j+1}^{0}\dots K_N^{0} v,  \label{e: e0}\\ 
& F_0 ( m \otimes v)    =   \sum_{j=1}^N m X_j^{-1} \otimes (K_{1}^{0})^{-1}\dots (K_{j-1}^{0})^{-1}E_j^{1,n} v,  \label{e: f0}\\  
&  K_0  = \kk_0 = (K_1 K_2 \cdots K_{n-1})^{-1}. \label{e: K0}
\end{align}
Here $ E_j^{i,k} = 1^{\otimes^{j-1}}\otimes E^{i,k} \otimes 1^{\otimes^{N-j}},$ where $E^{i,k} \in \End(V)$ is the matrix unit in the basis $ v_1,\dots,v_n,$ and $ K_j^i = q^{E_j^{i,i} - E_j^{i+1,i+1}},$ $ K_j^0 = (K_j^1 K_j^2 \cdots K_j^{n-1})^{-1}.$     

In view of the Theorem \ref{t: VV} and  the Proposition \ref{p: heckeact} the wedge product $\wedge^N V(z)$ is a left $\sltor$-module with $\kappa = q p^{-\frac{1}{n}}$ such that the action of the generators of the vertical and the horizontal subalgebras is given by the formulas (\ref{e: Efin} - \ref{e: K0}) where for  $m \in \cz$ and $v \in \otimes^N V$ we identify $m \otimes v$ with $\Lambda ( m \otimes v ),$ and use the maps in the  Proposition \ref{p: heckeact} to define a right action of $H_{tor}$ on $\cz.$ 

Sometimes we will denote the action of the vertical (horizontal) subalgebra on the wedge product $\wedge^N V(z)$ by $U_0^{(N)}$$(U_1^{(N)}).$

The actions of $E_{i,k}, \; F_{i,k}, \; H_{i,l}$ $(1 \leq i \leq n-1, \; k \in \zint , \; l \in \zint \setminus \{ 0 \} )$ are determined by the actions of the Chevalley generators of $\affaa $.

To prove the Theorem \ref{t: VV} 
 Varagnolo and Vasserot defined the following operator $\psi $ \cite{VV1}:
\begin{eqnarray}
& \psi (m \otimes v_{\ep _1} \otimes v_{\ep _2} \otimes \dots \otimes v_{\ep _N}) &  \label{psimap} \\ 
& =  m X_{1}^{-\delta _{n,\ep _1 }} X_{2}^{-\delta _{n,\ep _2}} \dots
X_{N}^{-\delta _{n,\ep _N}} 
 \otimes v_{\ep _1 +1}  \otimes v_{\ep _2 +1} \otimes \dots \otimes v_{\ep _N +1}. & \nonumber 
\end{eqnarray}
In this notation we identify $v_{n+1} = v_{1}$.
By taking into account the relations of the toroidal Hecke algebra, one can confirm that the action of $\psi $ is well-defined on $R\otimes_{\HH}(\otimes ^{N} V ).$

The following proposition is proved in \cite{VV1}.
\begin{prop} \label{vv3.4}
Let $E_{i}(z),F_{i}(z), K^{\pm } _{i}(z)$ be the series that specify actions of the generators of $\sltor$ on $R \otimes _{\HH} (\otimes ^{N} V )$. Then one has
\begin{eqnarray}
& \psi ^{-1} E_{i}(z) \psi = E_{i-1} (q^{-1} \kappa z),  
& \psi ^{-2} E_{1}(z) \psi ^{2}= E_{n-1} (x^{-1}q^{n-2} \kappa ^{2-n} z), \\
& \psi ^{-1} F_{i}(z) \psi = F_{i-1} (q^{-1} \kappa z),  
& \psi ^{-2} F_{1}(z) \psi ^{2}= F_{n-1} (x^{-1}q^{n-2} \kappa ^{2-n} z), \\
& \psi ^{-1} K^{\pm } _{i}(z) \psi = K^{\pm } _{i-1} (q^{-1} \kappa z),  
& \psi ^{-2} K^{\pm } _{1}(z) \psi ^{2}= K^{\pm } _{n-1} (x^{-1}q^{n-2} \kappa ^{2-n} z),
\end{eqnarray}
where $2 \leq i \leq n-1$.
\end{prop}

On the other hand, one can prove the following proposition \cite{VV1} easily.
\begin{prop} \label{psiprop}
Let $W$ be a $\affaa$-module such that actions of the generators are given by the series $E_i(z) , \; F_i(z) , \; H_i(z)$ $(1 \leq i \leq n-1)$  (\ref{defrel1} -- \ref{rg}). If there exists $\tilde{\psi } \in End(W)$ and $\kappa \in \cplx$ such that 
\begin{eqnarray}
& \tilde{\psi } ^{-1} E_{i}(z) \tilde{\psi } = E_{i-1} (q^{-1} \kappa z),  
& \tilde{\psi } ^{-2} E_{1}(z) \tilde{\psi } ^{2}= E_{n-1} (q^{-2} \kappa ^{2} z), \\
& \tilde{\psi } ^{-1} F_{i}(z) \tilde{\psi } = F_{i-1} (q^{-1} \kappa z),  
& \tilde{\psi } ^{-2} F_{1}(z) \tilde{\psi } ^{2}= F_{n-1} (q^{-2} \kappa ^{2} z), \\
& \tilde{\psi } ^{-1} K^{\pm } _{i}(z) \tilde{\psi } = K^{\pm } _{i-1} (q^{-1} \kappa z),  
& \tilde{\psi } ^{-2} K^{\pm } _{1}(z) \tilde{\psi } ^{2}= K^{\pm } _{n-1} (q^{-2} \kappa ^{2} z),
\end{eqnarray}
where $2 \leq i \leq n-1,$ 
then the $W$ is an  $\sltor $-module such that the actions of $E_0 (z), \; F_0 (z), \; H_0 (z)$ are given as follows:
\begin{eqnarray}
& E_0 (z) := \tilde{\psi } ^{-1} E_{1}(q \kappa ^{-1} z) \tilde{\psi }, & \\
& F_0 (z) := \tilde{\psi } ^{-1} F_{1}(q \kappa ^{-1} z) \tilde{\psi }, & \\
& K^{\pm } _0 (z) := \tilde{\psi } ^{-1} K^{\pm } _{1}(q \kappa ^{-1} z) \tilde{\psi }. & 
\end{eqnarray}
\end{prop}



\subsection{Semi-infinite wedge product}

In the subsection \ref{sec: qwedge}, we defined the space $V(z)$, its basis $\{ u _{k} \}, \; k \in \zint$ and the space $\wedge ^{N} V(z)$. 
In this section we define the semi-infinite wedge product $\wedge ^{\frac{\infty}{2}} V(z)$ and for any integer $M$ its subspace $F_M$ \cite{KMS}. Later we will define a representation of $\sltor$ on this space.

Let $\otimes ^{\frac{\infty}{2}} V(z)$ be the space spanned by the vectors
$u_{k_{1} } \otimes u_{k_{2}} \otimes \dots , \; ( k_{i+1} = k_{i} -1 , \; i >>1 )$. We define the space $\wedge ^{\frac{\infty}{2}} V(z)$ by the quotient of $\otimes ^{\frac{\infty}{2}} V(z)$:
\begin{equation}
\wedge ^{\frac{\infty}{2}} V(z) :=  \otimes ^{\frac{\infty}{2}} V(z) / \sum_{i=1}^{\infty }\Ker\left( \overset{c}T_i + q^2 (\overset{s}T_i)^{-1}\right). \label{e: inf-wedgeprod}
\end{equation}
Let $\Lambda: \otimes ^{\frac{\infty}{2}} V(z) \rightarrow \wedge ^{\frac{\infty}{2}} V(z)$ be the quotient map specified by (\ref{e: inf-wedgeprod}). The image of a pure tensor $ u_{k_1}\otimes u_{k_2}\otimes \dots  $ under this map is called a semi-infinite wedge and is denoted by
\begin{equation} 
u_{k_1}\wedge u_{k_2}\wedge \dots  :=\Lambda ( u_{k_1}\otimes u_{k_2}\otimes \dots ) . \label{e: infwedge} 
\end{equation}
A semi-infinite wedge is normally ordered if $k_1 > k_2 > \cdots $ and $k_{i+1}= k_i -1 \; (i >> 1)$. In \cite{KMS} it is proven that normally ordered semi-infinite wedges form a basis in $\wedge ^{\frac{\infty}{2}} V(z)$.
\vspace{.1in}

Let $U_{M}$ be the subspace of $\otimes ^{\frac{\infty}{2}} V(z)$ spanned by the vectors $u_{k_{1} } \otimes u_{k_{2}} \otimes \dots , \; ( k_{i} = M-i+1 , \; i >>1 )$.
Let $F_{M}$ be the quotient space of $U_{M}$ defined by the map (\ref{e: infwedge}).
 Then $F_{M}$ is a subspace of $\wedge ^{\frac{\infty}{2}} V(z)$, and the vectors $u_{k_1}\wedge u_{k_2}\wedge \dots $, ($k_1 > k_2 > \cdots $,  $k_{i} = M-i+1 , \; i >>1 )$ form the basis of $F_{M}$. We will  call the space $F_{M}$ the Fock space.

\section{The two actions of $U_q'(\sll)$ on the Fock space}
\subsection{Level-zero action of $U_q'(\sll)$ on the Fock space} \label{sec:lza}
Here we define a level-0 action of $U_q'(\sll)$ on $F_M$ $(M\in \zint)$ following the paper \cite{TU}. The definition we give below is equivalent to the one given in \cite{TU}. However, compared to \cite{TU}, we change slightly the precise wording so as to make the idea of this definition more transparent.   

Let $ \ee := (\ep_1,\ep_2,\dots,\ep_N) $ where $ \ep_i \in \{1,2,\dots,n\} $. For a sequence $\ee $ we set
\begin{equation}
{\bold v}_{\ee}:= v_{\ep_1}\otimes v_{\ep_2}\otimes \dots \otimes v_{\ep_N}\quad  ( \in \otimes^N \cplxn ).  
\end{equation}
A sequence $\mm := (m_1,m_2,\dots,m_N)$ from $\zint^N$ is called $n$-strict if it contains no more than $n$ equal elements of any given value. 
Let us define the sets $\MC_N^n$ and $\EC(\mm)$ by 
\begin{align}
 & \MC_N^n := \{\mm = (m_1,m_2,\dots,m_N) \in \zint^N \; | \; m_1\leq m_2 \leq \dots \leq m_N,\; \text{$\mm $  is $n$-strict } \}, \label{e: Mm}\\
 \intertext{ and for $ \mm \in \MC_N^n $ } 
& \EC(\mm) := \{ \ee = (\ep_1,\ep_2,\dots,\ep_N) \in \{1,2,\dots,n\}^N \; | \; \ep_i > \ep_{i+1} \; \text{for all $i$ s.t. $m_i =m_{i+1}$ } \}.   \label{e: Em}
\end{align}
In these notations the set 
\begin{equation}
\{ w(\mm,\ee):= \Lambda( \zz^{\mm} \otimes {\bold v}_{\ee}) \equiv z^{m_1}v_{\ep_1}\wedge z^{m_2}v_{\ep_2}\wedge\dots\wedge z^{m_N}v_{\ep_N} \quad | \quad \mm \in \MC_N^n , \ee \in \EC(\mm) \}.
\end{equation}
is nothing but the base of the normally ordered wedges in $\wedge^N V(z)$. We  will use the notation $w(\mm,\ee)$ {\em exclusively} for normally ordered wedges.

Similarly for a semi-infinite wedge $w = u_{k_1}\wedge u_{k_2} \wedge \dots \; $ $=$ $z^{m_1}v_{\ep_1}\wedge z^{m_2}v_{\ep_2} \wedge \dots \;,$ such that $w \in F_M,$ the semi-infinite sequences $\mm = (m_1,m_2,\dots \;\;)$ and $\ee = (\ep_1,\ep_2,\dots \;\;)$ are defined by $ k_i = \ep_i - nm_i $, $\ep_i \in \{1,2,\dots,n\},$ $m_i \in \zint.$ In particular the $\mm$- and $\ee$- sequences of the vacuum vector in $F_M$ will be denoted by $\mm^0$ and $\ee^0$:      
\begin{equation}
 |M\rangle = u_M\wedge u_{M-1} \wedge u_{M-2}\wedge  \dots \; = z^{m_1^0}v_{\ep_1^0}\wedge z^{m_2^0}v_{\ep_2^0} \wedge z^{m_3^0}v_{\ep_3^0} \wedge \cdots \;.  
\end{equation}

The Fock space $F_M$ is $\zintp$-graded. For any semi-infinite wedge $w$ $=$ $u_{k_1}\wedge u_{k_2} \wedge \dots \;$ $\equiv$  $z^{m_1}v_{\ep_1}\wedge z^{m_2}v_{\ep_2} \wedge \dots \;$ $\in F_M$ the degree $| w |$ is defined by  
\begin{equation}
| w | = \sum_{i \geq 1} m_i^0 - m_i . \label{wdeg}
\end{equation}
Later we will see that the $F_M$ is a toroidal module. In this module the degree generator $d_2$ acts as the grading operator whose eiegenvalue on a normally ordred wedge $w$ is equal to the degree (\ref{wdeg}). Clearly the degree is a finite non-negative integer for any wedge $u_{k_1}\wedge u_{k_2} \wedge \dots \;$ in $F_M$  because of the asymptotic condition $ k_i = M - i+ 1 ,$ $ i >> 1.$ Let us denote by $F_M^k$ $\subset$ $F_M$ the homogeneous component of degree $k$. 

We will define a level-zero action of $U_q'(\sll)$ on the Fock space $F_M$ in such a way that each homogeneous component $F_M^k$ will be invariant with respect to this action. Throughout this section we fix an integer $M$ and  $s \in \{0,1,2,\dots,n-1 \} $ such that  $M = s\bmod n.$ 

Let $l$ be a non-negative integer and define $ V_M^{s+nl} \subset \wedge^{s+nl}V(z) $ as follows:
\begin{equation}
V_M^{s+nl}  = \bigoplus\begin{Sb} \mm \in \MC_{s+nl}^n, \ee \in \EC(\mm) \\ m_{s+nl} \leq m_{s+nl}^0 \end{Sb} \cplx w(\mm,\ee) . \label{e: VM}
\end{equation}
Notice that the condition $m_{s+nl} \leq m_{s+nl}^0$ in this definition is equivalent to the condition
\begin{equation}
m_i \leq m_i^0 \quad \text{ for all $ i=1,2,\dots,s+nl$} \label{e: mlessm0}
\end{equation}
since the sequence $\mm$ is $n$-strict and non-decreasing.

The vector space $V_M^{s+nl}$ has a grading similar to the grading of the Fock space $F_M$. In this case the degree $|w |$ of a wedge  $w$ $=$ $u_{k_1}\wedge u_{k_2} \wedge \dots \wedge u_{k_{s+nl}}$ $\equiv$  $z^{m_1}v_{\ep_1}\wedge z^{m_2}v_{\ep_2} \wedge \dots \wedge z^{m_{s+nl}}v_{\ep_{s+nl}} $ $\in$  $V_M^{s+nl}$ is defined by     
\begin{equation}
| w | = \sum_{i = 1}^{s+nl} m_i^0 - m_i .
\end{equation}
Due to (\ref{e: mlessm0}) the degree is a non-negative integer, and for $k\in \zintp$ we denote by $ V_M^{s+nl,k} $ the homogeneous component of degree $k.$

The following result is contained in the paper \cite{TU}:
\begin{prop}
For each $k\in \zintp$ the homogeneous component $ V_M^{s+nl,k} $ $\subset$ $\wedge^{s+nl}V(z)$ is invariant under the $\affaa $-action  $U_0^{(s+nl)}$ defined in section $\ref{sec: finqta}$. 
\end{prop}

\begin{df}
For each $k \in \zintp$ define a map $\oro$ $:$ $V_M^{s+nl,k} $ $\rightarrow$ $F_M$ by setting for $w$ $\in$ $V_M^{s+nl,k}$ 
\begin{equation}
\oro(w)  =  w\wedge | M - s - nl \rangle .
\end{equation}
\end{df}
Clearly we have $|\oro(w)|$ $=$ $|w|$ and hence $\oro$ $:$ $V_M^{s+nl,k} $ $\rightarrow$ $F_M^k$ for all $k \in \zintp.$ 

\begin{prop}
When $l \geq k$ the map $\oro$ is an isomorphism of vector spaces. \label{p: isom}
\end{prop}
\begin{pf}
(Surjectivity) If a normally ordered wedge $w = u_{k_1}\wedge u_{k_2} \wedge \dots \; $ belongs to $F_M^k$ then $k_i = M-i+1$ for all $ i \geq s+nk+1.$ For if otherwise, the degree of $w$ must be greater or equal to $k+1$ (Cf. Proposition 5(ii) in \cite{TU}). Thus when $l \geq k$ we have $w$ = $w_{(s+nl)}\wedge |M-s-nl \rangle,$ and  $w_{(s+nl)} \in V_M^{s+nl,k}.$ Since a basis of $F_M^k$ is formed by normally ordered wedges, the surjectivity follows.

(Injectivity) If $w,\,w'$ $\in$ $V_M^{s+nl,k}$ are two distinct normally ordered wedges, then $w\wedge | M-s-nl \rangle,\, w'\wedge | M-s-nl \rangle $ are distinct and, as implied by the definition (\ref{e: VM}), normally ordered wedges in $F_M^k.$ Thus the injectivity follows. 
\end{pf} 

This proposition has an immediate corollary:
\begin{cor}
For each triple of non-negative integers $k,l,m$ such that $k\leq l < m$ the map $\orom$ $:$ $ V_M^{s+nl,k}$ $\rightarrow$ $ V_M^{s+nm,k},$ defined for any $w$ $\in$  $ V_M^{s+nl,k}$ by 
\begin{equation}
\orom(w) = w\wedge u_{M-s-nl}\wedge u_{M-s-nl-1}\wedge \dots \wedge u_{M-s-nm+1},
\end{equation}
is an isomorphism of vector spaces.  \label{c: cor1}
\end{cor}
Moreover we have a stronger statement:
\begin{prop} \label{u1isom}
For each triple of non-negative integers $k,l,m$ such that $k\leq l < m$ the map $\orom$ $:$ $ V_M^{s+nl,k}$ $\rightarrow$ $ V_M^{s+nm,k}$ is an isomorphism of the $U_q'(\sll)$-modules.  \label{p: isom1}
\end{prop}
\begin{pf}
In the proof of this and some of the subsequent propositions the following lemma concerning the Cherednik's operators and proved in the paper \cite{TU} plays an essential role:  
\begin{lemma}
Let $ \mm $ = $( m_1,m_2,\dots,m_N) \in \zint^N $ be a sequence such that:
\begin{equation} 
 m_1,m_2,\dots,m_{N-k} < m_{N-k+1} = m_{N-k+2}=\dots = m_N \; \equiv \; t ; \quad 1 \leq k \leq N.\end{equation}
Then  the following relations hold:
\begin{align}
&\text{\em{for} $0\leq l \leq k-1$} \qquad   \zz^{\mm} (Y_{N-l}^{(N)})^{\pm 1}  =  p^{\pm t} q^{\pm (2k-2l-N-1)} \zz^{\mm} + [\dots] \; , \label{e: ia}\\ 
&\text{\em{for} $1\leq i \leq N-k$} \qquad   \zz^{\mm} (Y_i^{(N)})^{\pm 1} =  q^{\pm k} \zz^{\mm} (Y_i^{(N-k)})^{\pm 1} + [\dots] \; , \label{e: ib}\\ 
\intertext{where $[\dots]$ signifies a linear combination of monomials $\zz^{\nn}$ $\equiv$ $ z_1^{n_1}z_2^{n_2}\dots z_N^{n_N}$ such that:}
& n_1,n_2,\dots,n_N  \leq  t ,  \qquad \text{and} \quad \#\{n_i | n_i = t\}  <  k.     
\end{align} \label{l: lemma3}
\end{lemma}

To prove the proposition it is sufficient to assume that $m$ is equal to $l+1$. And since the isomorphism of the linear spaces has been already established in the Corollary \ref{c: cor1} we must now show that for any generator $a$ of $U_q'(\sll)$ and any $w$ $\in$ $ V_M^{s+nl,k}$ the following intertwining relation takes place:  
\begin{equation}
a^{(s+nl+n)}\orol(w) = \orol a^{(s+nl)} (w) \qquad ( k \leq l ). \label{e: inter}
\end{equation}
Here $a^{(N)}$ represents the action of the generator $a$ in $U_0^{(N)}.$ 

Let $\mm' \in \MC_{s+nl}^n$ and $\ee' \in \EC(\mm')$ be such that $w(\mm',\ee')$ $\in$  $V_M^{s+nl,k}$  and take in (\ref{e: inter} ) $w$ $=$ $w(\mm',\ee')$ and $a = \ff _0.$  And let   $\mm \in \MC_{s+nl+n}^n$ and $\ee \in \EC(\mm)$ be such that
\begin{align}
& w(\mm,\ee) = \wedge(\zz^{\mm}\otimes {\bold v}_{\ee}) =  w(\mm',\ee')\wedge u_{M-s-nl}\wedge u_{M-s-nl-1}\wedge \dots \wedge u_{M-s-nl-n+1}, \\
& \text{ so that } \quad  w(\mm,\ee) = \orol(w(\mm',\ee')). 
\end{align}
Now act with  $\ff _0^{(s+nl+n)}$ on $w(\mm,\ee)$:
\begin{equation}
\ff _0^{(s+nl+n)}w(\mm,\ee) = \Lambda(\sum_{j=1}^{s+nl+n} q^{1-(s+nl+n)}\zz^{\mm} Y_j^{(s+nl+n)} \otimes
(K_1^0)^{-1}\dots(K_{j-1}^0)^{-1}E^{1,n}_j {\bold v}_{\ee}).
\end{equation}
Lemma \ref{l: lemma3} where we take $N=s+nl+n, k=n, t=m_{s+nl+n}^0$ and the normal ordering rules allow us to write 
\begin{multline}
\ff _0^{(s+nl+n)}w(\mm,\ee) = (\ff _0^{(s+nl)}w(\mm',\ee'))\wedge u_{M-s-nl}\wedge u_{M-s-nl-1}\wedge \dots \wedge u_{M-s-nl-n+1}  + \\ + p^{m_{s+nl+n}^0}q^{2(1-nl-n-s)}\,((\kk _0^{(s+nl)})^{-1}w(\mm',\ee'))\wedge \\ \wedge u_{M-s-nl-n+1}\wedge u_{M-s-nl-1}\wedge u_{M-s-nl-2}\wedge  \dots \wedge u_{M-s-nl-n+1} +  \tilde{w}, \label{e: lst}
\end{multline}
where the $\tilde{w}$ is a linear combination of normally ordered wedges $w(\nn,{\bold \tau}),$ $\nn \in \MC_{s+nl+n}^n,$ ${\bold \tau} \in \EC(\nn)$ such that for the sequence $\nn = (n_1,n_2,\dots,n_{s+nl+n})$ we have    
\begin{gather}
 n_1,n_2,\dots,n_{s+nl+n}  \leq  m_{s+nl+n}^0 ,   \\
\text{ and } \quad  \#\{n_i | n_i = m_{s+nl+n}^0 \}  <  n.     \label{e: less1}
\end{gather}
The inequality (\ref{e: less1}) implies, in particular, that $n_{s+nl+1} < m^0_{s+nl+1} = m^0_{s+nl+n}.$ And from the last inequality it follows that degree of the wedge  $w(\nn,{\bold \tau})$ is greater or equal to $l+1$ (Cf. Proposition 5(ii) in \cite{TU}). Since the degree of $ \ff _0^{(s+nl+n)}w(\mm,\ee)$ equals to the degree of $w(\mm,\ee)$ and equals to $k$, and since the degrees of the first two summands in (\ref{e: lst}) are equal to $k$ as well, taking the condition $k \leq l$ into account we find that $\tilde{w}$ equals to zero.

Now consider the second summand in (\ref{e: lst}). Lemma 2.2 in \cite{KMS} shows that
\begin{equation}
 u_{M-s-nl-n+1}\wedge u_{M-s-nl-1}\wedge u_{M-s-nl-2}\wedge  \dots \wedge u_{M-s-nl-n+1} = 0.
\end{equation}
Therefore we have 
\begin{equation}
\ff _0^{(s+nl+n)}w(\mm,\ee) = (\ff _0^{(s+nl)}w(\mm',\ee'))\wedge u_{M-s-nl}\wedge u_{M-s-nl-1}\wedge \dots \wedge u_{M-s-nl-n+1}, 
\end{equation}
which proves  (\ref{e: inter}) for $a$ $=$ $\ff _0.$

The proof of (\ref{e: inter}) for the rest of the $U_q'(\sll)$-generators is carried out in the same way.
\end{pf}

Now we make the final step and define on the vector space $F_M^k$ a level-zero action of $U_q'(\sll)$ by using Propositions \ref{p: isom} and \ref{p: isom1}:

\begin{df} 
The vector space $F_M^k$ is a level-$0$ module of $U_q'(\sll)$ with the action $U_0$ defined by  
\begin{equation}
 U_0 = \oro U_0^{(s+nl)} {\oro}^{-1} \quad \text{ where $ l \geq k .$ }  \label{eq:U0def}
\end{equation} 
\end{df}
Due to proposition \ref{p: isom1} this definition does not depend on the choice of $l$ as long as $l$ is greater or equal to $k.$ Since we have
\begin{equation}
F_M = \bigoplus_{k \geq 0} F_{M}^k
\end{equation}
the level-0 action $U_0$ extends to the entire Fock space $F_M.$

\subsection{Level-one action of $U_q'(\sll)$ on the Fock space} \label{KMSl1}

In this section we review the level-one action of $U_q'(\sll)$ on the Fock space $F_M$ \cite{KMS}. 

First we define the action of $U'_q(\sll)$ (generated by $E _{i}, \; F _i , \; K_i, \; $ $i=0, \dots ,n-1$) on the vector $|M' \rangle $ as follows.
\begin{equation}
E_i| M' \rangle =0,
\label{eq79}
\end{equation}
\vspace{-.18in}
\begin{equation}
F_i| M' \rangle = \left\{
\begin{array}{ll}
u_{M' +1}\wedge u_{M'-1}\wedge u_{M'-2}\wedge \cdots &  \mbox { if } i\equiv M' \mbox{ mod }n; \\
0& \mbox{ otherwise},
\end{array}
\right.
\end{equation}
\vspace{-.18in}
\begin{equation}
K_i | M' \rangle = \left\{
\begin{array}{ll}
q| M' \rangle & \mbox { if } i\equiv M' \mbox{ mod }n; \\
| M' \rangle & \mbox{ otherwise},
\end{array}
\right.
\end{equation}
For every element $v \in F_{M}$, there exist $N$ such that
\begin{equation}
v = v^{(N)} \wedge | M -N \rangle , \; \; \; \; \; 
v^{(N)} \in \wedge ^N V(z).
\end{equation}
We define the action of $E _{i}, \; F _i , \; K_i, \; $ $i=0, \dots ,n-1$ on the vector $v$ as follows.
\begin{eqnarray}
& E _i v := E _i v^{(N)} \wedge K_i | M -N \rangle +  v^{(N)} \wedge E_i | M -N \rangle ,& \\
& F _i v := F _i v^{(N)} \wedge | M -N \rangle +  K_i^{-1} v^{(N)} \wedge F_i | M -N \rangle ,& \\
& K _i v := K _i  v^{(N)} \wedge K _i | M -N \rangle.
\end{eqnarray}
The actions of $E _{i}, \; F _i , \; K_i, \; $ $i=0, \dots ,n-1$ on $v^{(N)}$ are determined in the section \ref{sec: finqta}. The definition of the actions on $v$ does not depend on $N$ and is well--defined, and we can easily check that the $U'_q (\sll)$-module defined in this section is level-1. We will use the notation $U_1$ for this  $U'_q (\sll)$-action on the Fock space.

\section{Action of $\sltor $ on the Fock space $F_M$}
\subsection{Action of $\sltor $ on the Fock space $F_M$}

On the module defined in \ref{sec: finqta}, the map (\ref{psimap}) 
can be written in the wedge notation as follows
\begin{equation}
\psi _{N} ( u_{k_1} \wedge u_{k_2} \wedge \dots \wedge u_{k_N} )
= u_{k_1 +1} \wedge u_{k_2 +1} \wedge \dots \wedge u_{k_N +1}
\end{equation}
For the toroidal  action on $\wedge^N V(z) = \otimes^N V(z)/ \sum_{i=1}^{N-1}\Ker\left( \overset{c}T_i + q^2 (\overset{s}T_i)^{-1}\right) $, we get the following relations
\begin{eqnarray}
& \psi _{N}^{-1} E_{i}(z) \psi _{N} = E_{i-1} (q^{-1} \kappa z),  
& \psi _{N}^{-2} E_{1}(z) \psi _{N}^{2}= E_{n-1} (q^{-2} \kappa ^{2} z), \\
& \psi _{N}^{-1} F_{i}(z) \psi _{N}= F_{i-1} (q^{-1} \kappa z),  
& \psi _{N}^{-2} F_{1}(z) \psi _{N}^{2}= F_{n-1} (q^{-2} \kappa ^{2} z), \\
& \psi _{N}^{-1} K^{\pm } _{i}(z) \psi _{N}= K^{\pm } _{i-1} (q^{-1} \kappa z),  
& \psi _{N}^{-2} K^{\pm } _{1}(z) \psi _{N}^{2}= K^{\pm } _{n-1} (q^{-2} \kappa ^{2} z),
\end{eqnarray}
 where $ i \in \{ 2,\dots,n-1\}$ and  $\kappa = p ^{-1/n}q$.

On the space $\wedge ^{\frac{\infty }{2}} V(z)$, we introduce the following map
\begin{equation}
\psi _{\infty} (u_{k_1} \wedge u_{k_2} \wedge \dots )
:= u_{k_1 +1} \wedge u_{k_2 +1} \wedge \dots \; . \label{psiact}
\end{equation}
Note that $\psi _{\infty } (F_M) = F_{M+1} $.

\begin{prop} \label{proppsiinf}
For each vector $v \in \wedge ^{\frac{\infty }{2}} V(z)$ and  $ i \in \{ 2,\dots,n-1\}$ we have 
\begin{eqnarray}
& \psi _{\infty }^{-1} E_{i}(z)\psi _{\infty }v = E_{i-1} (q^{-1} \kappa z)v 
,  
& \psi _{\infty }^{-2} E_{1}(z) \psi _{\infty }^{2}v= E_{n-1} (q^{-2} \kappa ^{2} z)v, \\
& \psi _{\infty }^{-1} F_{i}(z) \psi _{\infty }v= F_{i-1} (q^{-1} \kappa z)v,  
& \psi _{\infty }^{-2} F_{1}(z) \psi _{\infty }^{2}v= F_{n-1} (q^{-2} \kappa ^{2} z)v, \\
& \psi _{\infty }^{-1} K^{\pm } _{i}(z) \psi _{\infty }v= K^{\pm } _{i-1} (q^{-1} \kappa z)v,  
& \psi _{\infty }^{-2} K^{\pm } _{1}(z) \psi _{\infty }^{2}v= K^{\pm } _{n-1} (q^{-2} \kappa ^{2} z)v.
\end{eqnarray}
Here the action of $\affaa$ is the action $U_0$ defined in section \ref{sec:lza} and  $\kappa = p ^{-1/n}q$.
\end{prop}
\begin{pf}
To prove this proposition we need the following lemma.
\begin{lemma} \label{stab}
 Let $k,l,N$ be integers satisfying $l \geq k$ and $N=s+ln$. Here $s$ is the integer such that $ s \equiv M$ mod $n$, $0 \leq s \leq n-1$. Assume $v \in  V_M^{s+ln,k}$ then we have
\begin{eqnarray}
& E_i (z) (v \wedge u_{M-N} )= ( E_i (z) \cdot v ) \wedge u_{M-N} , & \label{E1} \\
& F_i (z) (v \wedge u_{M-N} )= ( F_i (z) \cdot v ) \wedge u_{M-N} ,& \\
& K^{\pm } _i (z) (v \wedge u_{M-N} )= ( K^{\pm } _i (z) \cdot v ) \wedge u_{M-N} ,& \label{K1}\\
& E_{n-1} (z) (v \wedge u_{M-N} \wedge u_{M-N-1})= ( E_{n-1} (z) \cdot v ) \wedge u_{M-N} \wedge u_{M-N-1}, & \label{En-1} \\
& F_{n-1} (z) (v \wedge u_{M-N} \wedge u_{M-N-1})= ( F_{n-1} (z) \cdot v ) \wedge u_{M-N} \wedge u_{M-N-1}, & \\
& K^{\pm } _{n-1} (z) (v \wedge u_{M-N} \wedge u_{M-N-1})= ( K^{\pm } _{n-1} (z) \cdot v ) \wedge u_{M-N} \wedge u_{M-N-1}, & \label{Kn-1}
\end{eqnarray}
where $1 \leq i \leq n-2$.
\end{lemma}
\begin{pf}
By the relations (\ref{rd}), (\ref{re}), we can show that for each $i$ ($1\leq i \leq n-1$) the subalgebra in $\affa$ generated by $E_{i,l'} , \; F_{i,l'}, \; H_{i,m'}, K^{\pm }_{i}$ $(l' \in \zint , \; m' \in \zint \setminus \{ 0 \} )$ is, in fact,  generated by only the elements $E_{i,0}, F_{i,0}, K^{\pm }_{i}, F_{i,1} $ and $F_{i,-1}$.
By the definition of the representation, every generator of $\affaa$ preserves the degree in the sense of (\ref{wdeg}).
So it is sufficient to show that the actions of $E_{i,0}, F_{i,0}, 
K^{\pm }_{i}, F_{i,1} $ and $F_{i,-1}$ satisfy the relations 
(\ref{E1}-- \ref{Kn-1}). For the case $E_{i,0}, F_{i,0}, K^{\pm }_{i}$,
by the definition of the actions (\ref{e: Efin}--\ref{e: Kfin}) we can check that that $E_{i,0}, F_{i,0}, K^{\pm }_{i}$ satisfy the relations (\ref{E1}-- \ref{Kn-1}) directly.
Let us show that 
\begin{eqnarray}
& F_{i,\pm 1} (v \wedge z^m v_n ) =  (F_{i,\pm 1} \cdot v) \wedge  z^m v_n, & \label{fi1} \\
& F_{n-1,\pm 1} (v \wedge z^m v_n \wedge z^m v_{n-1}) =  (F_{n-1,\pm 1} \cdot v) \wedge  z^m v_n \wedge z^m v_{n-1} , & \label{fn-11}
\end{eqnarray}
where $v \in  V_M^{s+ln,k}$  $(1 \leq i \leq n-2)$ and $m$ is such that $u_{M-(s+nl)} = z^m v_n$.

We will prove (\ref{fn-11}).
In the proof we will use the two different notations: 
\begin{align}
& u_{k_1}\wedge u_{k_2} \wedge \cdots \wedge u_{k_{N+2}} \quad \text{ or} \quad  \Lambda ( \zz^{\mm} \otimes {\bold v}_{\ee} ) \\ 
& \text{ where } \quad k_i = \ep_i - N m_i \quad  \text{ and } \quad \mm = (m_1,\dots,m_{N+2}), \quad \ee  = (\ep_1,\dots,\ep_{N+2})
\end{align}
for an element from $\wedge ^{N+2} V(z) = \cplx [ z_{1}^{\pm 1} , \cdots z_{N+2}^{\pm 1}] \otimes (\otimes ^{N+2} V) / \sum_{i=1}^{N+1} $Im$ (\stackrel{c}{T_i} - \stackrel{s}{T_i} )$.  



For any $ M',M'' $ ($1 \leq M' \leq M'' \leq N+2$), we define the $U'_q(\sll )$--action on the space $\cplx [z_{1}^{\pm 1}, \dots z_{N+2}^{\pm 1}] \otimes (\otimes ^{N+2} V)$ in terms of the Chevalley generators:
\begin{align}
& \ee_i( P(\zz ) \otimes w) =  \sum_{j=M'}^{M''} P(\zz ) \otimes E_j^{i,i+1} K_{j+1}^{i}\dots K_{M''}^{i} w, \label{eM} \\ 
& \ff_i ( P(\zz ) \otimes w)    =   \sum_{j=M'}^{M''}  P(\zz ) \otimes (K_{M'}^{i})^{-1}\dots (K_{j-1}^{i})^{-1}E_j^{i+1,i} w.  \\  
& \kk _i ( P(\zz ) \otimes w) =   P(\zz ) \otimes (K^i_{M'} K^i_{M'+1} \dots K^i_{M''}) w. \\
& \ee_0( P(\zz ) \otimes w) =  \sum_{j=M'}^{M''} P(\zz ) \cdot (q^{-N-1}Y_{j}^{(N+2)})^{-1} \otimes E_j^{n,1} K_{j+1}^{0}\dots K_{M''}^{0} w,  \\ 
& \ff_0 ( P(\zz ) \otimes w)    =   \sum_{j=M'}^{M''}  P(\zz ) \cdot (q^{-N-1} Y_{j}^{(N+2)}) \otimes (K_{M'}^{0})^{-1}\dots (K_{j-1}^{0})^{-1}E_j^{1,n} w.  \\  
& \kk _0 ( P(\zz ) \otimes w) =   P(\zz ) \otimes (K^0_{M'} K^0_{M'+1} \dots K^0_{M''}) w. \label{kM}   
\end{align}
Here $i= 1, \dots n-1$, $E^{l,l'}_{j} = 1 ^{\otimes ^{j-1}} \otimes E^{l,l'} \otimes 1 ^{\otimes ^{N+2-j}}, \; K^{i}_{j} = q^{E^{i,i}_{j} - E^{i+1,i+1}_{j}} $, $K^0_j = (K^1_j K^2_j \cdots K^{n-1}_j)^{-1}$, $P(\zz ) \in \cplx [z_{1}^{\pm 1}, \dots z_{N+2}^{\pm 1}]$ and $w \in  \otimes ^{N+2} V$.
This action is well--defined because of the commutativity of $Y_{i}^{(N+2)}$ ($i = 1, \dots , N+2$).
 The actions of the Drinfel'd generators are determined by the actions of the Chevalley generators.
 
Let $X$ be an element of $U'_q(\sll )$, then we define the action $X^{(*)}$ on the space $\cplx [z_{1}^{\pm 1}, \dots z_{N+2}^{\pm 1}] \otimes (\otimes ^{N+2} V)$ by (\ref{eM}--\ref{kM}) and $M'=1, \; M''=N$.

 We define the action $X^{(**)}$  on the space $\cplx [z_{1}^{\pm 1}, \dots z_{N+2}^{\pm 1}] \otimes (\otimes ^{N+2} V)$ by (\ref{eM}--\ref{kM}) and $M'=N+1, \; M''=N+2$.

 We also define the action $X^{\{ j \} }$ $(j=1, \dots , N+2)$ on the space $\cplx [z_{1}^{\pm 1}, \dots z_{N+2}^{\pm 1}] \otimes (\otimes ^{N+2} V)$ by (\ref{eM}--\ref{kM}) and $M'=j, \; M''=j$.

With these definitions, for any two elements $X$ and $Y$ from  $U'_q(\sll ),$  the operators  $X^{(*)}$ and $Y^{(**)}$ commute. This shows that if we have $\Delta X= \sum _{\nu } Y_{\nu } \otimes Z_{\nu }$ then $X (P(\zz ) \otimes w ) = \sum _{\nu } Y_{\nu }^{(*)} Z_{\nu }^{(**)} (P(\zz ) \otimes w ) $.

The following equations are satisfied modulo $ \Lambda (UN_+^{(*)} 
\cdot (UN_-^2)^{(**)} (P(\zz ) \otimes w ))$. Here $UN_+, \; 
UN_-^{2}$ are the left ideals generated by $\{ E_{i,k'} \} , \; 
\{ F_{i,k'}F_{j,l'} \}$.
\begin{align}
&  F_{n-1,1} \Lambda (P(\zz ) \otimes w ) \equiv \Lambda( (K_{n-1} ^{(*)}  F_{n-1,1}^{(**)}+  F_{n-1,1}^{(*)}) (P(\zz ) \otimes w_ )) , \\
&  F_{n-1,-1} \Lambda (P(\zz ) \otimes w ) \equiv \Lambda (((K_{n-1} ^{(*)})^{-1}  F_{n-1,-1}^{(**)}+  F_{n-1,-1}^{(*)}  \\
& \quad \quad \quad + (q^{-1}-q)(K_{n-1}^{(*)})^{-1} H_{n-1,-1}^{(*)}  F_{n-1,0}^{(**)} ) (P(\zz ) \otimes w )) ,\nonumber \\
& \quad  P(\zz ) \otimes w \in \cplx [z_{1}^{\pm 1}, \dots ,z_{N+2}^{\pm 1}] \otimes   (\otimes ^{N+2} V) . \nonumber
\end{align}
These equations follow from the the coproduct formulas which have been 
obtained in \cite{koyama} Prop. 3.2.A:
\begin{eqnarray}
& \Delta (F_{i,1}) \equiv K_{i} \otimes F_{i,1} + F_{i,1} \otimes 1 &
\mbox{mod } UN_+ \otimes UN_-^{2} , \label{cop1} \\
& \Delta (F_{i,-1}) \equiv  K_{i}^{-1} \otimes F_{i,-1} + F_{i,-1} \otimes 1 &
\label{cop2} \\
&  + (q^{-1}-q)K_i^{-1} H_{i,-1} \otimes F_{i,0} & \mbox{mod } UN_+
\otimes UN_-^{2}. \nonumber
\end{eqnarray}

Now we will show the equality 
\begin{equation}
\Lambda( F_{n-1,\pm 1} (P(\zz ) \otimes w )) 
=\Lambda(F_{n-1,\pm 1}^{(*)} (P(\zz ) \otimes w)).
\label{ffstar}
\end{equation}
where 
\begin{align}
 & P(\zz )= z_1^{m_1} z_2^{m_2} \cdots z_{N+2}^{m_{N+2}}, \; w = v_{\ep _1} \otimes \cdots  v_{\ep _{N+2}}, \text{ and } \label{Pvrel} \\
 & m_{N+1} = m_{N+2} =m, \; \; m_i <m \; ( i=1, \dots , N) , \nonumber \\ 
 & | \Lambda (P(\zz )  \otimes w))| = k , \quad \ep _{N+1}= n , \; \ep _{N+2}= n-1. \nonumber
\end{align}

First let us prove that any element in $UN_+ ^{(*)} \cdot (UN_-^2)^{(**)}$ 
annihilates  a vector $P(\zz ) \otimes w$ that  satisfies (\ref{Pvrel}).  
It is enough to show that 
\begin{equation}
 ( F_{i',k'}^{(**)}F_{j',l'}^{(**)})  
(z_{N+1}^m z_{N+2}^m \bar{v}  \otimes (v_n \otimes v_{n-1}))=0,
\end{equation}
for $\bar{v} \in \cz \otimes (\otimes ^N V)$.

This follows immediately from the observation that $wt(v_n)+wt(v_{n-1})-\alpha _{i'}-\alpha _{j'}$ is not a weight of $\otimes ^2 V$.

Next we will show that $\Lambda (F_{n-1,\pm 1}^{(**)} (P(\zz ) \otimes w) )= 0$, ($P(\zz )$ and $w$ satisfy (\ref{Pvrel})).
By the formulas (\ref{cop1}) and (\ref{cop2}), we have the following identities  modulo \\
$\Lambda (UN_+^{\{ N+1 \}  }  (UN_- ^2) ^{\{ N+2 \} } ( P (\zz ) \otimes w ))$:
\begin{align}
&  \Lambda ( F_{n-1,1}^{\{ N+2 \} } (P(\zz ) \otimes w )) \equiv \Lambda( (K_{n-1} ^{\{ N+1 \} }  F_{n-1,1}^{\{ N+2 \} }+  F_{n-1,1}^{\{ N+1 \} }) (P(\zz ) \otimes w )) , \\
&  F_{n-1,-1} \Lambda (P(\zz ) \otimes w ) \equiv \Lambda (((K_{n-1} ^{\{ N+1 \} })^{-1}  F_{n-1,-1}^{\{ N+2 \} }+  F_{n-1,-1}^{\{ N+1 \} }  \\
& \quad \quad \quad + (q^{-1}-q)(K_{n-1}^{\{ N+1 \} })^{-1} H_{n-1,-1}^{\{ N+1 \} }  F_{n-1,0}^{\{ N+2 \} } ) (P(\zz ) \otimes w )) ,\nonumber \\
& \equiv \Lambda (((K_{n-1} ^{\{ N+1 \} })^{-1}  F_{n-1,-1}^{\{ N+2 \} }+  F_{n-1,-1}^{\{ N+1 \} }  \nonumber \\
& \quad \quad \quad + (q^{-1}-q)[ E_{i,0}^{\{ N+1 \} }, F_{i,-1}^{\{ N+1 \} } ] F_{n-1,0}^{\{ N+2 \} } ) (P(\zz ) \otimes w )) ,\nonumber 
\end{align}
here we used the relation $[ E_{i,0}, F_{i,-1} ]= K_{i}^{-} H_{i,-1} $ which is proved by (\ref{re}). 
The following formula is essentially written in \cite{koyama} Prop. 3.2.B:
\begin{align}
& F_{i,\pm 1}^{\{ l\} } ( P'(\zz ) \otimes (\otimes _{j=1}^{N+2} v _{\ep _j}))  \label{ko2} \\
& =P'(\zz ) (q^{i}(q^{-N-1}Y_l ^{(N+2)})^{-1})^{\pm 1} \otimes (\otimes _{j=1}^{l-1} v _{\ep _j}) \otimes 
\delta _{i,\ep _l} v _{i+1} \otimes (\otimes _{j=l+1}^{N+2} v _{\ep _j}),
\nonumber
\end{align}
where $P'(\zz ) \in \cplx [ z_1^{\pm 1} , \dots ,z_{N+2}^{\pm 1}]$ and $v_j \in \cplx ^n$.

By the formula (\ref{ko2}), we have $(UN_+^{\{ N+1 \} } (UN_- ^2)^{\{ N+2 \} } (P(\zz ) \otimes w )) =0 $, and by the formula (\ref{ko2}), Lemma \ref{l: lemma3} (\ref{e: ia}) and the normal ordering rules, we get
\begin{align}
&  \Lambda (F_{n-1,\pm 1}^{(**)} (P(\zz ) \otimes w )  \label{ot} \\
& = \Lambda ( \alpha _{\pm 1} z_{N+1}^m z_{N+2}^m  \bar{v} \otimes v_n \otimes v_n) + \tilde{w } .\nonumber 
\end{align}
Here $ \alpha _{\pm 1}$ are certain coefficients, $ \bar{v} \in \cz \otimes ( \otimes ^N V)$ and $\tilde{w}$ is a linear combination of normally ordered wedges $w(\nn,{\bold \tau}),$ $\nn \in \MC_{s+nl+2}^n,$ ${\bold \tau} \in \EC(\nn)$ (see subsection \ref{sec:lza}) such that for the sequence $\nn = (n_1,n_2, \dots n_{s+nl+2})$ we have    
\begin{gather}
 n_1,n_2 , \dots ,n_{s+nl+2} \leq  m ,   
\text{ and } \quad  \#\{n_i | n_i = m \}  <  2. \label{ni2}   
\end{gather}
The inequality (\ref{ni2}) implies, in particular, that $n_{s+nl+1} < m.$ 
By the definition of the space $\wedge ^{N+2} V(z)$, we have $\Lambda ( z_{N+1}^m z_{N+2}^m  \bar{v} \otimes v_n \otimes v_n)=0 $.

Assume $\tilde{w}\neq 0$. From the  inequality  $n_{s+nl+1} < m$ it follows that 
$|\tilde{w}|\geq l+1$ (Cf. Proposition 5(ii) in \cite{TU}). On the other hand
we have $|\tilde{w}|=|\Lambda (P(\zz ) \otimes w)|=k$. By the condition
$k\leq l$ this is a contradiction. Therefore we conclude that $\tilde{w}=0,$ and hence $\Lambda (F_{n-1,\pm 1}^{(**)} (P(\zz ) \otimes w))=0$.

Let us prove $\Lambda((K^{(*)}_{n-1})^{-1} H^{(*)}_{n-1,-1} F_{n-1,0}^{(**)} 
(P(\zz ) \otimes w))=0.$

We have
\begin{align}
& \Lambda ((K^{(*)}_{n-1})^{-1} H^{(*)}_{n-1,-1} F_{n-1,0}^{(**)} (z_{N+1}^m z_{N+2}^m \bar{v} \otimes ( v_{n} \otimes v_{n-1}))) \\
& = \Lambda(( K^{(*)}_{n-1})^{-1} H^{(*)}_{n-1,-1} ( z_{N+1}^m z_{N+2}^m \bar{v} \otimes ( v_n \otimes v_{n}))) , \nonumber
\end{align}
here $ \bar{v} \in \cz \otimes ( \otimes ^N V)$.
Since $H^{(*)} _{n-1,-1}$ belongs to the algebra generated by the operators $\ee _i ^{(*)}, \; \ff ^{(*)}_i , \; (\kk ^{(*)}_i)^{\pm } $ $(i=0, \dots , n-1)$, $H_{n-1,-1}^{(*)}$ belongs to the algebra generated by the operators  $(Y_{j}^{(N+2)})^{\pm 1}$,  $E_j ^{l,l'}  \; (1 \leq j \leq N , \; 1 \leq l,l' \leq n)$. 
By Lemma \ref{l: lemma3} (\ref{e: ib}) and the normal ordering rules, we have
\begin{align}
& \Lambda( (K^{(*)}_{n-1})^{-1} H^{(*)}_{n-1,-1} (z_{N+1}^m  z_{N+2}^m \bar{v} \otimes ( v_n \otimes v_{n}))) \\
& = \Lambda( z_{N+1}^m  z_{N+2}^m \tilde{v} \otimes ( v_n \otimes v_{n})) + \tilde{w }. \nonumber
\end{align}
Here $\tilde{v} \in \cz \otimes (\otimes ^N V)$ and $\tilde{w}$ is the element which has the property written after the relation (\ref{ot}). By the previous discussion we have $\tilde{w}=0$.

Thus we have shown (\ref{ffstar}).
 To prove (\ref{fn-11}) we must show that in the right-hand side of the last 
equation we can replace $q^{1-(N+2)}Y^{(N+2)}_i$ by $q^{1-N}Y^{(N)}_i$ $(1 \leq i \leq N).$ By  Lemma \ref{l: lemma3} (\ref{e: ib}), for $1\leq i \leq N$ and a sequence $\mm = (m_1 , m_2, \dots m_{N+2}) \in \zint ^{N}$ such that $m_1 , m_2 , \dots m_N < m_{N+1} = m_{N+2}=m $ we have:
\begin{equation}
\zz^{\mm} (Y_i^{(N+2)})^{\pm 1}  =  q^{\pm 2} \zz^{\mm} (Y_i^{(N)})^{\pm 1} + [\dots] ,
\end{equation}
where $[\dots]$ signifies a linear combination of monomials $\zz^{\nn}$ $\equiv$ $ z_1^{n_1}z_2^{n_2}\dots z_{N+2}^{n_{N+2}}$ such that $n_1,n_2,\dots,n_{N+2}  \leq  m $ and $\#\{n_i | n_i = m \}  <  2$. By the normal ordering rules, we can write
\begin{align}
& \Lambda (q^{1-(N+2)}Y^{(N+2)}_i(z_{N+1}^{m} z_{N+2}^{m} \bar{v} \otimes (v_n \otimes v_{n-1}))) \\
& = \Lambda (q^{1-N}Y^{(N)}_i (z_{N+1}^{m} z_{N+2}^{m} \bar{v} \otimes (v_n \otimes v_{n-1}))) + \tilde{w}, \nonumber
\end{align}
where  $\bar{v} \in \cz \otimes (\otimes ^N V)$ and the $\tilde{w }$ again has the same meaning as the $\tilde{w }$ in  relation (\ref{ot}).  Repeating the discussion after (\ref{ot}) we can show that $\tilde{w}=0$.
Hence we get (\ref{fn-11}).
\vspace{.2in}

To prove (\ref{fi1}), consider the tensor product $\cplx [z_1 ^{\pm 1} , \cdots , z_{N+1} ^{\pm 1}] \otimes (\otimes ^{N} V) \otimes V$, use the formulas (\ref{cop1}), (\ref{cop2}) and continue  the  proof in a way that is  completely analogous to the proof of (\ref{fn-11}).
\end{pf}
\vspace{.2in}

Now we proceed with the proof of Proposition \ref{proppsiinf}. 
It is sufficient to show the statement of the proposition for the vector $v$ such that $v \in F^{k}_{M}$. 
We put $v= v_{(N)} \wedge | M-N \rangle $, where $v_{(N)} \in V_{M}^{N,k}, \; N=s+ln $ as in Lemma \ref{stab}.
By  Proposition \ref{psiprop}, Lemma \ref{stab} (\ref{E1}) - (\ref{K1}) and the relation (\ref{psiact}), we can show the following relations ($X= E, \; F $ or $K^{\pm 1} $, $2 \leq i \leq n-1$):  

\begin{align}
& X_{i-1}(q^{-1} \kappa z) \cdot v \label{Xv}\\
& = ( X_{i-1} (q^{-1} \kappa z) v_{(N)} ) \wedge  u_{M-N} \wedge |M-N-1 \rangle \nonumber \\
& =  X_{i-1} (q^{-1} \kappa z) (v_{(N)} \wedge  u_{M-N} ) \wedge |M-N-1 \rangle  \nonumber \\
& = ( \psi _{N+1}^{-1} X_{i}(z) \psi _{N+1} (v_{(N)} \wedge  u_{M-N} )) 
\wedge |M-N-1 \rangle  \nonumber \\
& = \psi _{\infty }^{-1}  X_{i}(z) (\psi _{N+1} (v_{(N)} \wedge  u_{M-N} ) \wedge |M-N \rangle  )  \nonumber \\
& = \psi _{\infty }^{-1}  X_{i}(z) (\psi _{\infty } ( v_{(N)} \wedge  u_{M-N}  \wedge |M-N-1 \rangle  ),  \nonumber
\end{align}
This is exactly the statement of the Proposition \ref{proppsiinf} ($i \neq n-1$). For the case $i=n-1$ we use Lemma \ref{stab} (\ref{En-1}) - (\ref{Kn-1}) and a straightforward modification of (\ref{Xv}).
\end{pf}
\vspace{.2in}

In view of  Propositions \ref{psiprop} and \ref{proppsiinf} we can now define a $\sltor $-action on the space $\wedge ^{\frac{\infty}{2}} V(z)$. By this definition  each subspace $F_M$ is invariant with respect to the $\sltor $-action. Thus we arrive at  
\begin{thm}[Varagnolo and Vasserot \cite{VV2}] \label{mthm}
The q-deformed Fock space $F_M$ is a $\sltor $-module.
The actions of $X_i(z) \; (X=E,F$ or $K^{\pm }, \; 1 \leq i \leq n-1)$ are determined in subsection \ref{sec:lza} by the Chevalley generators. The actions of $X_0(z) \; (X=E,F$ or $K^{\pm })$ are determined by $ \psi _{\infty }^{-1} X_1(p^{1/n}z)  \psi _{\infty }$ (in this notation $X_1(p^{1/n}z)$ act on $F_{M+1}$, but $X_0(z)$ act on $F_M$).
\end{thm}
\begin{rmk}
One can verify that the action of the subalgebra $\affbb $ coincides with the level-1 action of $U_q'(\sll)$ defined in Section \ref{KMSl1}, because we have $X_{0} = \psi _{\infty }^{-1} X_1  \psi _{\infty }$ $(X=E,F$ or $K^{\pm })$ ($X_i$ is the Chevalley generator of the level-1 action of $U_q'(\sll)$ defined in Section \ref{KMSl1}). Hence the action of  $\sltor $ on the Fock space has the level $(0,1)$.
\end{rmk}

\subsection{Action of $\sltor $ on level-1 irreducible $U_q'(\sll)$-modules}

In the paper \cite{KMS} it was demonstrated that the Fock space $F_M$ admits an action of the Heisenberg algebra $H$ which commutes with the level-1 action  $U_1$ of the algebra $U_q'(\sll).$ The Heisenberg algebra is a unital $\cplx$-algebra generated by elements $1, B_a$ with $a \in \zint_{\neq 0}$ which are subject to  relations 
\begin{equation}
[ B_a , B_b ] = \delta_{a+b,0}a\frac{1 - q^{2n a}}{1 - q^{2 a}} . \label{e: Heisenbergrel} \end{equation}
The Fock space $F_M$ is an $H$-module with the action of the generators given by \cite{KMS} 
\begin{equation}
B_a = \sum_{i=1}^{\infty} z_i^{a}. \label{e: Bact}
\end{equation}

Let $\cplx[H_{-}]$ be the Fock space of $H$, i.e.,  $\cplx[H_{-}]= \cplx[B_{-1},B_{-2},\dots,\:].$ The element $B_{-a}$ $ ( a = 1,2,\dots\:)$ acts on  $\cplx[H_{-}]$ by multiplication. The action of  $B_{a}$ $ ( a = 1,2,\dots\:)$ is given by (\ref{e: Heisenbergrel}) together with the relation
\begin{equation}
 B_a \cdot 1 = 0 \qquad \text{ for $ a \geq 1. $} 
\end{equation}

Let $\Lambda_i$ $(i\in \{ 0,1,\dots,n-1\}) $ be the fundamental weights of $\sll'.$ And let $V(\Lambda_i)$ be the irreducible (level-1) highest weight module of $U_q'(\sll)$ with highest weight vector $V(\Lambda_i)$ and highest weight $\Lambda_i.$

The following results are proven in \cite{KMS}:
\begin{itemize}
\item The action of the Heisenberg algebra on $F_M$ and the action $U_1$ of $U'_q (\sll )$ commute. 

\item There is an isomorphism 
\begin{equation}
 \iota_M : F_M \cong V(\Lambda_i)\otimes \cplx[H_{-}] \qquad ( M= i\bmod n )
\label{e: iota} \end{equation}
of $U_q'(\sll)\otimes H$-modules normalized so that $\iota_M( |M\rangle )$ $ =$ $V(\Lambda_i)\otimes 1.$
\end{itemize}
In general the level-0 $U_q'(\sll)$-action $U_0$   does not commute with the Heisenberg algebra. However if we choose the parameter $p$ in $U_0$ in a special way, then the following result holds:   
\begin{prop}
At $p=1$ we have
\begin{equation}
[ U_0 , H_{-}] = 0.
\end{equation}  \label{p: U0comH-}
\end{prop}
\begin{pf}
Let $w \in F_M^k.$ Then by Proposition \ref{p: isom} for any $l \geq k$ there is a unique $w_{(s+nl)} \in V_M^{s+nl,k}$ such that 
\begin{equation}
w = w_{(s+nl)} \wedge | M - s - nl \rangle . \label{e: w=w}
\end{equation}
Let $m \geq 1$ and let us act with $B_{-m}$ on the $w$: 
\begin{equation}
B_{-m} w = ( B_{-m}^{(s+nl)} w_{(s+nl)})\wedge | M - s - nl \rangle +  w_{(s+nl)}\wedge B_{-m} | M - s - nl \rangle ,  \label{e: Bw}
\end{equation}
where $ B_{-m}^{(N)} = \sum_{i=1}^N z_i^{-m} $ $( N =1,2, \dots \:).$ In view of (\ref{e: w=w}) and (\ref{e: Bact}) in the second summand above we may write 
\begin{align}
& w_{(s+nl)} = w_{(s+nk)}\wedge u_{M-s-nk}\wedge u_{M-s-nk-1}\wedge \dots u_{M-s-nl+1}, \\
& B_{-m} | M - s - nl \rangle = u_{M-s-nl + nm}\wedge u_{M-s-nl-1} \wedge  u_{M-s-nl-2}\wedge \dots \;  + \\
& + u_{M-s-nl}\wedge u_{M-s-nl+nm-1} \wedge  u_{M-s-nl-2}\wedge \dots \; + \nonumber \\ 
& + u_{M-s-nl}\wedge u_{M-s-nl-1} \wedge  u_{M-s-nl+nm-2}\wedge \dots \; + \nonumber \\ 
& + \cdots \;. \nonumber
\end{align}
Now let us choose $l$ to be greater or equal to $ k + m.$ Then by Lemma 2.2 in \cite{KMS} the second summand in (\ref{e: Bw}) vanishes and  we have
\begin{equation}
B_{-m} w = ( B_{-m}^{(s+nl)} w_{(s+nl)})\wedge | M - s - nl \rangle.  \label{e: Bw1}
\end{equation}
The degree of $B_{-m} w$ is equal to $k+m.$ Note that since $m \geq 1$ we have also  
\begin{equation}
B_{-m}^{(s+nl)} w_{(s+nl)} \in V_M^{s+nl,k+m}. \label{e: In}
\end{equation}
Therefore by definition (\ref{eq:U0def}) of the action $U_0$ for any element $a$ of $U_q'(\sll)$ we have (note that $l \geq k + m$): 
\begin{equation}
a B_{-m} w = \left( a^{(s+nl)} B_{-m}^{(s+nl)} w_{(s+nl)} \right)\wedge | M - s - nl \rangle.
\end{equation}
At $p=1$ the operators $ a^{(N)}$ and $ B_{a}^{(N)}$ commute for all finite $N$ and non-zero $a$ since  $ B_{a}^{(N)}$ are symmetric in $z_1,\dots,z_N$ and thereby commute with the operators $Y^{(N)}_1,\dots,Y^{(N)}_N.$ 

Thus we have 
\begin{equation}
a B_{-m} w = \left( B_{-m}^{(s+nl)} a^{(s+nl)} w_{(s+nl)} \right)\wedge | M - s - nl \rangle.
\end{equation}
Taking into account that degree of $a w$ is equal to the degree of $w$ and is equal to $k$ we may repeat the discussion leading to (\ref{e: Bw1}) and find that  
\begin{equation}
a B_{-m} w =  B_{-m} a w.
\end{equation}
\end{pf}

\begin{rmk}
Note that it is not true that $U_0$ at $p=1$ commutes with the subalgebra $H_{+}$ generated by $1, B_1,B_2,\dots\;.$ This fact is known from consideration of the Yangian limit in \cite{U}. If we attempt to repeat the discussion in the proof of the preceding proposition for $B_{-m}$ with negative $m$, we find that the inclusion  (\ref{e: In}) fails to hold in general.    
\end{rmk}

Let $H_{-}'$ be the non-unital subalgebra in $H$ generated by  $ B_{-1},B_{-2},\dots\;.$ Proposition \ref{p: U0comH-} allows us to define a level-0 $U_q'(\sll)$-module structure on the irreducible level-1 module $V(\Lambda_i)$ $(i \in \{0,1,\dots,n-1\}).$ Indeed from this proposition it follows that the subspace  
\begin{equation}
  H_{-}'F_M \subset F_M 
\end{equation}
is invariant with respect to the action $U_0$ at $p=1$ and therefore a level-0 action of $U_q'(\sll)$ is defined on the quotient space
\begin{equation}
F_M/ (H_{-}'F_M)  \label{eq: jjj}
\end{equation}
which in view of (\ref{e: iota}) is isomorphic to  $V(\Lambda_i)$ with $ i = M\bmod n.$ We do not know whether this level-0 action coincides with the level-0 action defined in the paper \cite{JKKMP}. However the Yangian limit considered in \cite{U} suggests an affirmative answer to this question. The results of \cite{U}  also lead to the following conjecture: 
\begin{con}
At $p= q^{2n}$ we have
\begin{equation}
 [ U_0 , H_{+} ] = 0.
\end{equation}
\end{con}

According to Lemma \ref{l: gen} the subalgebras $U_1$ and $U_0$ generate the action of $\sltor$ on the Fock space, hence combining Proposition \ref{p: U0comH-} with the fact that $U_1$ commutes with the Heisenberg algebra, we find that the toroidal action at $p=1( \kappa = q )$ commutes with $H_{-}.$ Repeating the discussion leading to the equation  (\ref{eq: jjj}) we conclude that   
\begin{prop}
The highest weight irreducible $U_q'(\sll)$ module $V(\Lambda_i)$ $(i=0,\dots,n-1)$ is a level $(0,1)$ toroidal module with $\kappa = q.$ In this module the action of the subalgebra $\affbb$ coincides with the level $1$ action of $U_q'(\sll).$ And the action of the  subalgebra $\affaa$ is induced from the level $0$ $U_q'(\sll)$-action $U_0$ on the Fock space.   
\end{prop}

\subsection{Irreducibility of the Fock space as a toroidal module}

At $q=1$ both $\affaa $ and $\affbb $ are isomorphic to $U'(\sll).$ Action of these subalgebras  on the Fock space $F_M$ is now given by the generators  
\begin{gather}
E_i  =  \ee_i    = \sum_{j=1}^{\infty}  E_j^{i,i+1} , \qquad   
F_i    =  \ff_i    = \sum_{j=1}^{\infty} E_j^{i+1,i} ,  \label{e: Ffinc}\\  
H_i    =  \hh_i   =  \sum_{j=1}^{\infty} {E_j^{i,i} - E_j^{i+1,i+1}}   ,\qquad    (i=1,2,\dots,n-1)  \label{e: Kfinc} \\ 
\ee _0 =  \sum_{j=1}^{\infty}p^{-D_j} E_j^{n,1} ,  \qquad 
\ff _0    =   \sum_{j=1}^{\infty} p^{D_j} E_j^{1,n}, 	 \label{e: F0c} \\
E_0 =  \sum_{j=1}^{\infty}z_j  E_j^{n,1} , \qquad 
F_0  =   \sum_{j=1}^{\infty} z_j^{-1}E_j^{1,n} ,  \label{e: f0c}\\  
\hh_0 =  - \hh_1 - \hh_2 - \cdots - \hh_{n-1}, \label{e: h0c} \\
H_0  = 1 - H_1 - H_2 - \cdots - H_{n-1}. \label{e: H0c}
\end{gather}
It is straightforward to verify that these generators are well-defined  on the Fock space provided we specify normalization of the Cartan subalgebra in $\sln$ as  
\begin{equation}
H_i | M \rangle     =  \hh_i | M \rangle  = \delta( M \equiv i\bmod n )| M \rangle \qquad  (i=1,2,\dots,n-1). 
\end{equation}
The affine Lie algebra $\sll'$ is realized as the central extension of the loop algebra $ \sln\otimes\cplx[t,t^{-1}]$ by the center $\cplx c,$ so that if $ x(m) = x\otimes t^m$ for $x \in \sln,$ then a system of generators for  $U'(\sll)$ is provided by $ E_{i,j}(m)$ $ (i\neq j \in \{1,\dots,n\}),$ $ H_i(m)$ $  (i \in \{1,\dots,n-1\},\: m \in \zint)$  and $c,$ where $E_{i,j}$ are matrix units regarded as generators of $ \sln$ and $ H_i = E_{i,i} -E_{i+1,i+1}.$ 

In terms of these generators the action of $\affbb$ is given by 
\begin{gather}
 E_{i,j}(m) =  \sum_{k=1}^{\infty}z_k^m  E_k^{i,j} , \quad   
H_i(m)  =   \sum_{k=1}^{\infty}z_k^m ( E_k^{i,i} -  E_k^{i+1,i+1}) ,  \quad c=1. \label{e: gen}  
\end{gather}
\begin{prop}
Let $q=1$ and $|p| \neq 1.$ Then the Fock space $F_M$ is an irreducible module of the algebra  $\sltor.$ 
\end{prop}
\begin{pf} 
Let $B_m $ $( m\in \zint\setminus\{0\})$ be the generators (\ref{e: Bact}) of the $H$-action on $F_M.$ A computation shows that we have   
\begin{equation}
B_m = \sum_{i=1}^{n-1}(i + \frac{p^m}{1-p^m} ) H_i(m) + \frac{n p^{m}}{p^{2m} -1} [[ E_{1,n}(0), [ E_{n,1}(m) , \ff _0 ]], \ee _0 ] 
\end{equation}
and therefore the Heisenberg action is included into the action of toroidal algebra. This implies that the $\affbb \otimes H$ action is included into the toroidal action as well. However the former action is irreducible in view of the decomposition   (\ref{e: iota}).  
\end{pf}
\begin{cor}
Let $|p| \neq 1$. Then there is $\ep > 0$ such that  the Fock space $F_M$ is an irreducible module of the algebra  $\sltor$  for all $ q $ with $  |q - 1| < \ep.$
\end{cor}
\vspace{.2in}

{\bf Acknowledgement}
We thank Professors Masaki Kashiwara and Tetsuji Miwa for support and discussions.


\begin{thebibliography}{99}

\bibitem{BPS}
Bernard, D., Pasquier, V. and Serban D.: Spinons in conformal field theory,  Nucl. Phys., {\bf B428} 612-628 (1994).
\bibitem{CP} 
Chari, V. and Pressley, A.:  Quantum affine algebras and affine Hecke
 algebras, Pacific Jour. of Math. {\bf 174} 295-326 (1996) (q-alg/9501003).
\bibitem{Cherednik}
Cherednik, I.V.: A unification of the Knizhnik-Zamolodchikov and Dunkl operators via affine Hecke algebras, Inv. Math., {\bf 106} 411-432 (1991). 
\bibitem{Cherednik2} 
Cherednik, I.V.: Double Affine Hecke algebras and Macdonald's conjectures, Annals Math., {\bf 141} 191-216 (1995); Non-symmetric Macdonald's Polynomials, IMRN {\bf 10} 483-515  (1995).
\bibitem{gkv}
Ginzburg, V., Kapranov, M. and Vasserot, E.:
Langlands reciprocity for algebraic surfaces,
Math. Res. Lett. {\bf 2} 147-160 (1995).
\bibitem{H}
Hayashi, T.: Q-analogues of Clifford and Weyl algebras -- spinor and oscillator representations of quantum enveloping algebras,  Comm. Math. Phys., {\bf 127} 129-144 (1990).
\bibitem{JKKMP}
Jimbo, M., Kedem, R., Konno, H., Miwa, T. and Petersen, J.-U. H.: Level-0 action of $U_q(\slt)$ on level-1 modules and Macdonald Polynomials,  J.Phys, {\bf A28} 5589 (1995) (q-alg/9506016).
\bibitem{KMS}
Kashiwara, M., Miwa, T. and Stern, E.: Decomposition of $q$-deformed Fock Spaces, Selecta Mathematica, New Series,  {\bf 1} No. 4 , 787-805 (1995). (q-alg/9508006)
\bibitem{koyama}
Koyama, Y.:  Staggered polarization of vertex models with $U_q(\sll )$-symmetry, Comm. Math. Phys., {\bf 164} 277-291 (1994).
\bibitem{saito}
Saito, Y.: Quantum toroidal algebras and their vertex representations,
Preprint RIMS-1112 (q-alg/9611030).
\bibitem{TU} Takemura, K. and Uglov, D.: Level-0 action of $U_q(\sll )$ on the q-deformed Fock spaces, Preprint RIMS-1096 (q-alg/9607031).
\bibitem{U}
Uglov, D.: Semi-infinite wedges and the conformal limit of the fermionic spin Calogero-Sutherland model of spin $\frac{1}{2},$ Nucl. Phys. {\bf B478} 401-430 (1996).
\bibitem{VV1} Varagnolo, M. and Vasserot, E.: Schur duality in the toroidal setting, Comm. Math. Phys., {\bf 182} 469-484 (1996).
\bibitem{VV2} Varagnolo, M. and Vasserot, E.: Double-loop algebras and the Fock space,  Preprint (q-alg/9612035).
\end{thebibliography}
\end{document}